\documentclass[11pt]{article}
\usepackage{a4}
\usepackage{graphicx}
\usepackage{parskip}

\usepackage{amsfonts}
\usepackage{natbib}
\usepackage{algorithm}
\usepackage{algorithmic}

\usepackage{latexsym}

\newtheorem{theorem}{Theorem}[section]
\newtheorem{proposition}{Proposition}[section]
\newtheorem{lemma}{Lemma}[section]
\newtheorem{corollary}{Corollary}[section]
\newtheorem{remark}{Remark}[section]
\newtheorem{example}{Example}[section]
\newtheorem{definition}{Definition}[section]

\newcommand{\vvert} {| \hskip-0.15em |  \hskip-0.15em |}
\newcommand{\R}{\mathbb{R}}
\newcommand{\Nat}{\mathbb{N}}

\newcommand{\PP} {{  \rm I\hskip-0.22em P}}

\newcommand{\EE} {{\rm I\hskip-0.48em E}}

\usepackage{a4}

\usepackage{graphicx}

\usepackage{parskip}

\usepackage{amsfonts}

\usepackage{natbib}

\usepackage{algorithm}

\usepackage{algorithmic}

\usepackage{amsmath, amssymb}

\usepackage{latexsym}

\begin{document}

\centerline{\bf \Large On the asymptotic variance of the debiased Lasso}

\vskip .1in
\centerline{Sara van de Geer}

\centerline{Seminar for Statistics, ETH Z\"urich}

\vskip .1in
\centerline{August 21, 2018} 

{\bf MSC 2010 subject classifications:} Primary 62J07; secondary 62E20.\\
{\bf Keywords and phrases:} asymptotic efficiency, asymptotic variance, Cram\'er Rao lower bound, 
debiasing, Lasso, sparsity. 

{\bf Abstract.} We consider the high-dimensional linear regression model
$Y = X \beta^0 + \epsilon$
with Gaussian noise $\epsilon$ and Gaussian random design $X$.
We assume that $\Sigma:= \EE X^T X / n$ is non-singular and write its inverse
as $\Theta := \Sigma^{-1}$. The parameter of interest is
the first component $\beta_1^0$ of $\beta^0$. 
We show that in the high-dimensional case the asymptotic variance of a debiased Lasso
estimator can be smaller than $\Theta_{1,1}$.  For some special such cases we establish
asymptotic efficiency. The conditions include
$\beta^0$ being  sparse 
and the first column $\Theta_1$ of $\Theta$ being not sparse.
These conditions depend on whether $\Sigma$ is known or not.

\section{Introduction} \label{introduction.section}

Let $Y$ be an $n$-vector of observations and $X \in \R^{n \times p} $ an input matrix.
The linear model is
$$ Y = X \beta^0 + \epsilon , $$
where $\beta^0 \in \R^p$ is a vector of unknown coefficients and $\epsilon \in \R^n$
is unobservable noise.  
We examine the high-dimensional case with
$ p \gg n$. The parameter of interest in this paper is a component of $\beta^0$, say the first 
component $\beta_1^0$. We consider the asymptotic properties of debiased
estimators of the one-dimensional
parameter $\beta_1^0$ under scenarios
where certain sparsity assumptions fail to hold. 

The paper shows that the asymptotic variance of the debiased estimator
can be smaller than the ``usual" value for the low-dimensional case. For simplicity we examine Gaussian
data: the rows of $(X,Y) \in \R^{n\times (p+1)}$ are i.i.d.\ copies
of a zero-mean Gaussian row vector $({\bf x}, {\bf y} ) \in \R^{p+1}$, where ${\bf x}= ({\bf x}_1 , \ldots ,
{\bf x}_p)$ has covariance matrix $\Sigma := \EE {\bf x}^T {\bf x} $. We assume the inverse
of $\Sigma$ exists and write it as
$\Theta:= \Sigma^{-1}$. The vector $\beta^0$ of regression coefficients is 
$ \beta^0 = \Theta \EE {\bf x}^T {\bf y} $.
We denote the first column of $\Theta$ by  $\Theta_1 \in \R^p$ and the first
element of this vector by $\Theta_{1,1}$. Our main aim is to present
examples where lack of sparsity in $\Theta_1$ can result in a small asymptotic
variance of a suitably debiased estimator. In particular, the asymptotic
variance can be smaller than $\Theta_{1,1}$. For the case $\Sigma$ known, this means applying for instance
a (noiseless) node-wise Lasso, instead of an exact orthogonalization of the first variable with respect to the others,
may reduce the asymptotic variance (as follows from combining Theorem \ref{bsharp.theorem} with Lemma \ref{findpair.lemma}).
If $\Sigma$ is unknown, the high dimensionality of the problem
already excludes exact empirical projections for orthogonalization. The (noisy) Lasso is designed 
to deal with approximate orthogonalization in the high-dimensional case.
Using the node-wise Lasso, we find that
one may again profit from non-sparsity of the (now unknown) vector $\Theta_1$
(see Theorem \ref{noisy.theorem}).

We look at specific examples or constructions of
covariance matrices $\Sigma$. 
The results illustrate that 
  asymptotic efficiency claims require some caution.
  The high-dimensional situation exhibits new phenomena
that do not occur in low dimensions. 

Throughout, the minimal
eigenvalue of $\Sigma$, denoted by $\Lambda_{\rm min}^2$, is required to
stay away from zero, i.e.\ $1/\Lambda_{\rm min}^2 = {\mathcal O} (1)$.
We further consider only $\Sigma$'s with all 1's on the diagonal and
assume for simplicity that $\sigma_{\epsilon}^2 := \EE
( {\bf y} - {\bf x} \beta^0 )^2$ is known and that its value is $\sigma_{\epsilon}^2=1$. 
   
Let a given subset
${\cal B}$ of $\R^p$ be
the model class for $\beta^0$. An interesting research goal is to construct for the model ${\cal B}$ a regular 
estimator of $\beta_1^0$ with asymptotic variance that achieves the asymptotic Cram\'er Rao lower bound
(given here in Proposition \ref{CRLB.proposition}). One then needs to decide which model class ${\cal B}$ 
one considers as relevant.
In high-dimensional
statistics it is commonly assumed that $\beta^0$ is sparse in some sense. 
Let $0< r \le 1 $, define  for a vector $b \in \R^p$ its $\ell_r$-``norm"
$\| b \|_r$ by $\| b \|_r^r := \sum_{j=1}^p | b_j |^r $ and let $\| b \|_0^0 $ be its number
of non-zero entries.
A sparse model  is for example
\begin{equation} \label{sparse-model.equation}
{\cal B} := \{ b \in \R^p : \ \| b \|_0^0 \le s \} 
 \end{equation}
 for some (``small") $s \in \Nat$.
 Alternatively one may believe only in $\ell_1$-sparsity. Then
 \begin{equation} \label{l1-model.equation}
 {\cal B} = \{ b \in \R^p : \ \| b \|_1 \le \sqrt {s} \}
 \end{equation}
 for some $s >0$.
These are the two extremes of weakly sparse models of the form
\begin{equation} \label{lr-model.equation}
 {\cal B} = \{ b \in \R^p : \| b \|_r^r \le s^{2-r \over 2}\} , 
 \end{equation}
 for some $s>0$ and $0\le r \le1 $. Throughout, the value of $s$ is allowed to depend on $n$, but
 $r$ is fixed for all $n$. 
  
 Constructing estimators that achieve the asymptotic Cram\'er Rao lower bound for model (\ref{sparse-model.equation}), ({\ref{l1-model.equation}), (\ref{lr-model.equation}) or some other sparse model
 is to our understanding
 quite ambitious, especially if one wants to do this for all possible covariance matrices $\Sigma$.
 See e.g.\ Example \ref{sparsemodel.example} for some details concerning model ({\ref{sparse-model.equation}).
 However, for special cases of $\Sigma$'s the problem can be solved.
  One such special case is where the first row of $\Theta_1$ of $\Theta$ is sparse in an appropriate sense. 
 This is the situation considered in previous work such as  \cite{zhang2014confidence} and \cite{van2013asymptotically}
 where $\Sigma$ is unknown. In this paper we consider known
 and unknown $\Sigma$ and in both cases do not require sparsity of $\Theta_1$. 
 The paper \cite{javanmard2018debias} also does not require sparsity of $\Theta_1$
 when $\Sigma$ is known and 
it turns out that for certain non-sparse vectors $\Theta_1$ their estimator
is not asymptotically efficient (see Theorem \ref{bsharp.theorem} or Remark \ref{efficient.remark}
following this theorem).


 The debiased Lasso defined in this paper in equation (\ref{de-biased.equation}) below is based on
  a direction $\tilde \Theta_1 \in \R^p$ where typically $\tilde \Theta_1$ is thought of as some estimate of $\Theta_1$. 
 As we do not
 assume sparsity of $\Theta_1$ a reliable estimator of $\Theta_1$ may not be available.
 Nevertheless, we show that this does not rule out good theoretical performance. 
 We present a class of covariance matrices $\Sigma$ for which a debiased
 Lasso has asymptotic variance smaller than $\Theta_{1,1}$. 
 This phenomenon is tied to the high-dimensional situation, see Remark \ref{high-dim.remark}.
 For special cases,
we establish that
 an asymptotic Cram\'er Rao bound smaller than $\Theta_{1,1}$ can  be achieved.
 In other words,  there exist cases where a
 debiased Lasso profits from sparsity of $\beta^0$ combined with non-sparsity of
$\Theta_1$. This is good news: the asymptotic variance can be small for two reasons:
either $\Theta_1$ is sparse in which case the asymptotic variance $\Theta_{1,1}$ is typically small
(as it is the inverse of the residuals of regressing the first variable  on
only a few of the other variables), or $\Theta_1$ is not sparse
 but then the asymptotic variance can be smaller than $\Theta_{1,1}$. This paper presents
 cases where the latter situation indeed occurs.  
   
 \subsection{The Lasso, debiased Lasso and sparsity assumptions}\label{Lasso.section}

The Lasso (\cite{tibs96}) is
defined as 
\begin{eqnarray}\label{Lasso.equation}
\hat \beta := \arg \min_{b \in \R^p} \biggl \{ \| Y - X b \|_2^2/n + 2 \lambda \| b \|_1 \biggr \}  
\end{eqnarray}
with $\lambda>0$ a tuning parameter. (We will throughout take of order $\sqrt {\log p / n }$ but not too
small.)

A debiased Lasso is
given by
\begin{equation}\label{de-biased.equation}
\hat b_1 = \hat \beta_1 + \tilde \Theta_1^T X^T (Y- X \hat \beta)/n .
\end{equation}
The $p$-dimensional vector $\tilde \Theta_1 $ is typically some estimate of the first column 
$\Theta_1$ of 
the precision matrix $\Theta$, but in our case it will rather be estimating a sparse approximation. 
We refer to $\tilde \Theta_1$ as a direction. The estimator $\hat \beta$ is commonly taken to be the Lasso given
in (\ref{Lasso.equation}) although this is not a must.

 The debiased Lasso (\ref{de-biased.equation}) was introduced in
\cite{zhang2014confidence} and further developed in \cite{javanmard2014confidence} and 
\cite{van2013asymptotically} for example.  Related work is \cite{belloniuniform}  and \cite{belloni2016post}.

The theory for the Lasso (\ref{Lasso.equation}) and debiased Lasso (\ref{de-biased.equation}) 
requires some form of sparsity of $\beta^0$.
Consider for some $s \in \Nat$ one of the sparsity models (\ref{sparse-model.equation}),
(\ref{l1-model.equation}) or, more generally, model (\ref{lr-model.equation}). 
 Two prevalent sparsity assumptions are:\\
(i) $s = o ( n / \log p )$,\\
(ii) $s = o ( \sqrt n / \log p )$.\\

 Sparsity variant (i) is for example invoked to establish $\ell_2$-consistency of the Lasso $\hat \beta$
(see \cite{bickel2009simultaneous} or the monographs \cite{koltchinskii2011oracle}, \cite{BvdG2011} and \cite{giraud2014introduction},
and their references). 

Which sparsity variant
is needed to establish asymptotic normality of the debiased Lasso  (\ref{de-biased.equation}) depends to a large extent
on whether $\Sigma$ is known or not. In \cite{javanmard2014hypothesis}, \cite{cai2017confidence}, \cite{ren2015}. 
\cite{javanmard2018debias} one can find
refined results on this issue.

This case  $\Sigma$ known is treated in Section \ref{knownSigma.section}. We introduce and apply there the concept of an eligible pair, see Definition \ref{gammasharp.condition}. An eligible pair is a sparse approximation
of $\Theta_1$ together with a parameter describing the order of approximation and sparsity. 
We allow for sparsity variant (i) in model (\ref{sparse-model.equation}) as in \cite{javanmard2018debias}, see Example
\ref{l0.example}. Sparsity variant (i) will also be allowed for the models
(\ref{l1-model.equation}) and (\ref{lr-model.equation}), see Examples \ref{l1.example} and
\ref{lr.example}.

 Eligible pairs will also play a crucial role in Section \ref{unknownSigma.section}
 where $\Sigma$ is unknown.  Let us discuss some of the literature for this case
 and for the sparsity model
(\ref{sparse-model.equation}).
 From the papers
\cite{cai2017confidence} and \cite{ren2015}  we know that for the minimax bias of an
estimator of $\beta_1^0$ to be of order $1/ \sqrt n$, the assumption $s=  O( \sqrt {n / \log p} )$ is necessary.
 Thus, up to log-terms  this needs 
 the second sparsity  variant.
 When considering asymptotic Cram\'er Rao lower bounds,
 one also needs to restrict oneself to a certain class of estimators,
 for instance estimators with bias of small order $1/\sqrt n$ or asymptotically linear 
 estimators. In \cite{Jankova16} such restrictions are studied.
  One can show asymptotic linearity of the debiased Lasso under model (\ref{sparse-model.equation}) with sparsity variant
  (ii)  and in addition $\| \Theta_1 \|_0 = o (\sqrt n / \log p )$. 
   If $\Theta_1$ is not sparse
  nor can be approximated by a sparse vector, then 
 it is unclear whether an asymptotically linear estimator exists. 
 We refer to Remark \ref{asymptotic-linearity.remark} for more details.
 In summary, modulo log-terms, sparsity variant (ii) cannot be relaxed as far as minimax rates
 for the bias are concerned, and sparsity variant (ii) with in addition sparsity
 of order $o (\sqrt n / \log p ) $  for $\Theta_1$ or its sparse approximation 
appears to be needed
for establishing asymptotic linearity.
We note that  the paper \cite{javanmard2018debias}
  establishes asymptotic normality under (among others) the assumption
  \begin{equation}\label{montenari.equation}
 \min ( s, \| \Theta_1 \|_0^0) = o( \sqrt n / \log p ).
 \end{equation}
 Bias and asymptotic linearity
 are not considered (these issues are not within the scope of that paper).
 In our setup however, $\Theta_1$ is not sparse at all, so
 variant (ii) is in line with (\ref{montenari.equation}). 

Tables \ref{conditionsknown.table} and \ref{conditionsunknown.table} presented in Subsection \ref{organization.section}
summarizes the sparsity conditions applied in this paper. One sees that models
(\ref{sparse-model.equation}) and (\ref{l1-model.equation}) are special cases of model
(\ref{lr-model.equation}), with $r=0$ and $r=1$ respectively. However, when $r=0$ the
asymptotic efficiency depends on $\beta^0$ and also quite severely on the value of $s$.
For the case $\Sigma$ unknown, model (\ref{l1-model.equation}) is too large. 

\subsection{The asymptotic Cram\'er Rao lower bound}\label{CRLB.section}
We briefly review the Cram\'er Rao lower bound and refer to \cite{Jankova16} for details.
Let the model be $\beta^0 \in {\cal B}$, where ${\cal B}$ is a given class of
regression coefficients. Let  ${\cal H}_{\beta^0} := \{ h \in \R^p: \  \beta^0 + h/\sqrt n   \in {\cal B} \} $.
We call ${\cal H}_{\beta^0} $ the set of model directions. 
An estimator $T$ (or actually: sequence of estimators) is called regular at $\beta^0$ if for all fixed $\rho>0$
and $R >0$ and all sequences $h  \in {\cal H}_{\beta^0}$ with $| h_1 | \ge \rho$ and
$h^T \Sigma h \le R^2$,
it holds that
$$ \sqrt n \biggl ( {T- (\beta_1^0+ h_1/\sqrt n)  \over  V_{\beta^0} } \biggr ){\buildrel{{\cal D}_{ \beta^0 + h/\sqrt n   }}\over \longrightarrow} {\cal N} (0,1)  $$
where $V_{\beta^0}^2 = {\mathcal O} (1)$ is some constant (depending on $n$ and possibly on $\beta^0$, but not depending
on $\rho$, $R$ or $h$), called the asymptotic variance
(it is defined up to smaller order terms). Regularity is important in practice. It means that the
asymptotics is not just pointwise but remains valid in neighbourhoods.

\begin{remark} Typically, the class ${\cal B}$ is not a cone. This is the reason why
we do not restrict ourselves to model directions $h \in {\cal H}_{\beta^0}$ with $h_1=1$ (say).
\end{remark}

\begin{remark} One may also opt for defining the set of possible directions 
${\cal H}_{\beta^0}$ differently,
say  ${\cal H}_{\beta^0} := {\cal B}$.
Then regularity concerns parameter values that fall outside the parameter space. 
For example under model (\ref{sparse-model.equation}) one then has to
deal with sparsity $2s$ instead of $s$. With our choice of ${\cal H}_{\beta^0}$ we stay inside the parameter space
but the set of possible directions then depends on $\beta^0$. In model (\ref{l1-model.equation}) or
(\ref{lr-model.equation}) one can move away from this dependence when $\beta^0$ is
a proper ``interior point" of ${\cal B}$ (see Examples \ref{l1-model.example} and
\ref{lr-model.example}). 
\end{remark}

\begin{proposition}\label{CRLB.proposition}
Suppose $T$ is asymptotically linear at $\beta^0$ with influence function
${\bf i}_{\beta^0}: \R^{p+1} \rightarrow \R$:
$$ T- \beta_1^0 = {1 \over n} \sum_{i=1}^n {\bf i}_{\beta^0} (X_i , Y_i) + o_{\PP_{\beta^0 }} (1/ \sqrt n ) $$
where $\EE_{\beta^0}{\bf i}_{\beta^0}({\bf x} , {\bf y} ) = 0 $ and
$V_{\beta^0}^2 := \EE_{\beta^0} {\bf  i}_{\beta^0}^{2} ({\bf x} , {\bf y} )= {\mathcal O} (1)$.
Assume the Lindeberg condition
$$\lim_{n \rightarrow \infty} \EE_{\beta^0} {\bf i}_{\beta^0}^2 ({\bf x} , {\bf y} ) {\rm l} \biggl \{ {\bf i}_{\beta^0}^2 ({\bf x} , {\bf y} )> \eta {n V_{\beta^0}^2} \biggr \} =0  \ \forall \ \eta >0 .  $$
Assume further that $T$ is regular at $\beta^0$. Then for all fixed $\rho>0$ and $R>0$
$$V_{\beta^0}^2 + o(1) \ge  \max_{h \in {\cal H}_{\beta^0}: \ |h_1| \ge \rho,\ 
h^T \Sigma h \le R^2} { h_1^2 \over
h^T \Sigma h } . $$
\end{proposition}

This proposition is as Theorem 9 in  \cite{Jankova16} but tailored for the particular situation.
A proof is given in Section \ref{proofs.section}.
We remark that such results are not a direct consequence of the Le Cam theory, as
we are dealing with triangular arrays. 

\begin{corollary}\label{CRLB.corollary} Assume the conditions of Proposition \ref{CRLB.proposition} and
that for some fixed $\rho>0$ and $R >0$ and some sequence $h \in {\cal H}_{\beta^0}$, with
$|h_1| \ge \rho$ and $h^T \Sigma h \le R^2$,  it is true that
$$ V_{\beta^0}^2 = { h_1^2 \over h^T \Sigma h } + o(1) . $$
Then $T$ is asymptotically efficient. This is our main tool to arrive at asymptotic efficiency
for some special $\Sigma$'s. 
\end{corollary}

The restriction to directions in ${\cal H}_{\beta^0}$ means that the Cram\'er Rao lower bound for
the asymptotic variance $V_{\beta^0}^2$ can be orders of magnitude smaller than $\Theta_{1,1}$.

\begin{example} \label{sparsemodel.example} Under the sparse model (\ref{sparse-model.equation})\footnote{A more natural model in this context might be ${\cal B} := \{ b \in \R^p : \ \| b_{-1} \|_0^0 \le s-1 \}$ where,
for $b \in \R^p$, $b_{-1} := ( b_2 , \ldots , b_p)^T$.}
 we have
$$ \underline {\cal H}_{\beta^0} \subseteq {\cal H}_{\beta^0} , $$
where
$$  \underline {\cal H}_{\beta^0}  := \{ \| h \|_0^0 \le s- s_0 \} $$
and $s_0 = |S_0|$ with 
$S_0 := \{  \beta_j^0 \not= 0 \}  $ being the set of active coefficients of $\beta^0$.

{\bf Some special cases}\\
a) If $\Theta_1  \in \underline {\cal H}_{\beta^0}$, then $V_{\beta^0}^2 + o(1)= \Theta_{1,1} $.
Note that the condition on $\Theta_1$ depends on $\beta^0$ (via $s_0$).\\
b)  Suppose $\{ 1 \} \in S_0$, $s=s_0$ and that the following ``betamin" condition
holds: $| \beta_j^0 | > m_n/\sqrt n$ for all $j \in S_0$, where $m_n $ is some sequence
satisfying 
$m_n \rightarrow \infty$. Then we see that 
$${\cal H}_{\beta^0} \cap \{ h^T \Sigma h \le R^2  \} =\{  \| h_{-S_0 } \|_0 = 0 \} \cap \{ h^T \Sigma h \le R^2  \}.$$
The lower bound is then
$$V_{\beta^0}^2+ o(1) \ge ( \Sigma_{S_0  , S_0 }^{-1})_{1,1} , $$
where $\Sigma_{S_0 , S_0 }$ is the matrix of covariances of the variables in $S_0 $.
This lower bound corresponds to the case $S_0$ is known.
The bound could be achieved if one has a consistent estimate $\hat S$ of $S_0$.
For this one needs betamin conditions in order to have no false negatives.
Estimation of $\beta_1^0$ using some estimator $\hat S$ of $S_0$ typically means that the estimator 
of $\beta_1^0$ is not regular.
Imposing further conditions beyond model (\ref{sparse-model.equation})
to make it regular
may diminish the lower bound.\\
c) More generally, if  $\{ 1 \} \in S_0$ and $| \beta_j^0 | > m_n /\sqrt n$ for all $j \in S_0$ and some
sequence $m_n \rightarrow \infty$ then the lower bound corresponds to knowing
the set $S_0$ up to $s-s_0$ additional variables. \\
d) Suppose that $\beta^0$ is an interior point in the sense that it stays away from the boundary:
for some fixed $0 < \eta< 1$ it holds that $s_0 \le (1- \eta) s $
(so that $1- s_0 / s$ stays away from zero). Then 
$$ \{ h : \ \| h \|_0^0 \le \eta s \} \subset {\cal H}_{\beta^0 }  .$$
The lower bound can then still depend on $\eta$.
A rescaling argument as applied in the next examples, Example \ref{l1-model.example} and
more generally Example \ref{lr-model.example},
 does not work here.
This signifies that
model (\ref{sparse-model.equation}) less suitable in our context: the exact value of $s$ plays a too
prominent role. 
\end{example}

\begin{example}\label{l1-model.example}
Consider model (\ref{l1-model.equation}). The situation is then more like a classical one.
Suppose $\beta^0$ stays away from the boundary of the parameter space, i.e., for a fixed $0<\eta < 1$
it holds that
$\| \beta^0 \|_1 \le (1-\eta )\sqrt s$. Then
$$ \{ \| h \|_1 \le \eta \sqrt {ns} \} \subset {\cal H}_{\beta^0 }  .$$
By a rescaling argument the dependence on $\eta$ in the left hand side plays no role in the lower bound:
for all $M>0$ fixed
$$ V_{\beta^0}^2 +o(1) \ge \biggl ( \min_{c \in \R^{p-1}: \ \| c \|_1 \le M\sqrt {ns} }\EE ( {\bf x}_1 - {\bf x}_{-1} c )^2 \biggr )^{-1} , $$
where ${\bf x}_1$ is the first entry of ${\bf x}$ and ${\bf x}_{-1} := ({\bf x}_2 , \ldots , {\bf x}_p )$. 
In other words, unlike in model (\ref{sparse-model.equation}), the exact value of $s$ does not play a prominent role.
We thus see that in order to be able to improve over $\Theta_{1,1}$ we must have that
$\Theta_1$ is rather non-sparse: $\| \Theta_1 \|_1$ should be of larger order
than $\sqrt {ns} $. 
For example if
$s = o(n /\log p)$, say $s= n/(m_n \log p )$ for some sequence $m_n \rightarrow \infty$ slowly,
the model directions have $\ell_1$-norm of order $n / \sqrt {m_n \log p} $.
To improve over $\Theta_{1,1}$ one thus must have $\|\Theta_1\|_1$ of larger
order, say $\|\Theta_1\|_1= n / \sqrt {\log p } $. We get back to this
in Example \ref{l1.example}. 

\end{example} 

\begin{example} \label{lr-model.example}
We now turn to model (\ref{lr-model.equation}). Using the same arguments as
in Example \ref{l1-model.example} one sees that if $\beta^0$ stays away from the boundary,
one may use model directions  with $\| \cdot \|_r^r $ of order $\sqrt {n^r s^{2-r}}$.
The lower bound is then
$$ V_{\beta^0}^2 +o(1) \ge \biggl ( \min_{c \in \R^{p-1}: \ \| c \|_r^r \le M\sqrt {n^r s^{2-r}} }\EE ( {\bf x}_1 - {\bf x}_{-1} c )^2 \biggr )^{-1}  $$
where $M >0$ is any fixed constant. 

\end{example}

It is clear, and illustrated by Examples \ref{sparsemodel.example} , \ref{l1-model.example} and
\ref{lr-model.example}, that the
lower bound of Proposition \ref{CRLB.proposition} depends on the model ${\cal B}$. 
The sparse model (\ref{sparse-model.equation}) is perhaps too
stringent. One may want to take the model ${\cal B}$ as the largest set for which a regular estimators exist.
This points in the direction  of  model (\ref{l1-model.equation}). We 
will see that when $\Sigma$ is known this model is indeed useful but when $\Sigma$ is unknown it is
too large.

\subsection{Notations and definitions}\label{notation.section}

We consider an asymptotic framework with triangular arrays of observations.
Thus, unless explicitly stated otherwise, all quantities depend on $n$ although we do not (always) express this in our notation.

Let ${\bf x}_1$ be the first entry of ${\bf x}$ and ${\bf x}_{-1}:= ({\bf x}_2 , \ldots, {\bf x}_p)$ be this vector with the first entry excluded,
so that ${\bf x} = ( {\bf x}_1 , {\bf x}_{-1} )$. Write $\Sigma_{-1,-1}:= \EE {\bf x}_{-1}^T {\bf x}_{-1} \in
\R^{(p-1)\times (p-1)} $. 
For vectors $b \in \R^p$ be use a similar
notation: $b_1\in \R$ is the first coefficient and $b_{-1} \in \R^{p-1}$ forms the rest of
the coefficients.
Apart from the regression (projection) ${\bf x} \beta^0$ of ${\bf y}$ on ${\bf x}$, we consider the
  regression ${\bf x}_{-1} \gamma^0$ of ${\bf x}_1$ on ${\bf x}_{-1}$, and an approximation
  $\gamma^{\sharp}$ of $\gamma^0$ which is accompanied by a parameter $ \lambda^{\sharp}$
  to form an ``eligible pair" $(\gamma^{\sharp}, \lambda^{\sharp})$ as given in Definition \ref{gammasharp.condition}.
  When $\Sigma$ is known we can invoke the noiseless Lasso $\gamma_{\rm Lasso}$ with tuning parameter
  $\lambda_{\rm Lasso}$ 
  (an approximate projection of ${\bf x}_1$ on
  ${\bf x}_{-1}$) to approximate $\gamma^0$. See (\ref{noiselessLasso.equation}) for its definition.
    For the case $\Sigma$ is unknown we apply the notation $\hat \Sigma := X^T X / n$. We let
  $X_1\in \R^n$ be the first column of $X$, and
the matrix $X_{-1}\in \R^{n \times (p-1)} $ be the remaining columns and we write $\hat \Sigma_{-1,-1}:=
X_{-1}^T X_{-1} / n $. 
We
do a approximate regression of $X_{1}$ on $X_{-1}$
 invoking the noisy Lasso $\hat \gamma$ with tuning parameter
  $\lambda^{\rm Lasso}$ as given in (\ref{noisyLasso.equation}). The various vectors of coefficients
  and their ``lambda parameter" are summarized in Table \ref{parameters.table}.
  Here we also add  the Lasso $\hat \beta$  for the estimation of $\beta^0$, as
  defined in (\ref{Lasso.equation}) with
  tuning parameter $\lambda$. 
      
  \begin{table}[!htbp]
\centering
\begin{tabular}{c|clclclcl}
 \ & \ & lambda     \\
 \ & coefficients& parameter \\
\hline
projection ${\bf y}$ on ${\bf x} $ & $\beta^0 $ & \ \ 0 \\
\hline
noisy Lasso & $\hat \beta$ & \ \ $\lambda \asymp \sqrt {\log p / n} $ \\
\hline
projection ${\bf x}_1$ on ${\bf x}_{-1}$ & $\gamma^0 $ & \ \ 0 \\
\hline
eligible pair & $\gamma^{\sharp} $& \ \ $\lambda^{\sharp}$ \\
\hline
noiseless Lasso & $\gamma_{\rm Lasso} $& \ \ $ \lambda_{\rm Lasso}$  \\
\hline
noisy Lasso & $\hat \gamma$ & \ \ $ \lambda^{\rm Lasso}\asymp \sqrt {\log p / n} $ \\
  \end{tabular}
   \caption{The various coefficients and lambda parameters}
   \label{parameters.table}
\end{table}

We write for $S \in \{1 , \ldots , p \}$, ${\bf x}_S :=\{ {\bf x}_j\}_{j \in S} $ and
${\bf x}_{-S} :=\{ {\bf x}_j\}_{j \notin S, \ j \not= 1} $. 
We further let for $S \subset \{ 1 , \ldots , p\}$
with cardinality ${\rm s}$, the matrix $\Sigma_{S,S}:= \EE {\bf x}_S^T {\bf x}_S \in \R^{{\rm s} \times {\rm s} }$ be the 
covariance sub-matrix formed
by the variables in $S$ and
$\Sigma_{-S,-S}:= \EE {\bf x}_{-S}^T {\bf x}_{-S} \in \R^{(p-1-{\rm s}) \times( p-1-{\rm s})}$ and 
$\Sigma_{S,-S} := \EE {\bf x}_S^T{\bf x}_{-S}  =: \Sigma_{-S, S}^T \in \R^{{\rm s} \times (p-1-{\rm s})}$.

Note that the vector of coefficients $\gamma^0$ of the
regression of ${\bf x}_{1}$ on ${\bf x}_{-1}$ is a vector in $\R^{p-1}$.
We will index its entries by $\{ 2 , \ldots , p \}$: $\gamma^0 = (\gamma_2^0 , \ldots , \gamma_p^0 )^T $.
Also the various other ``gamma" parameters" $\gamma$ will be indexed by $\{ 2, \ldots , p \}$.
It should be clear from the context when this indexing applies. 
For $c=(c_2, \ldots , c_p)$
and $S \subset \{ 2 , \ldots , p \}$ we write
$c_S := \{\gamma_j : \ j \in S \} $ and $\gamma_{-S} := \{ c_j :\ j \notin S, \  j \not= 1  \}$.
We use the same notation for the $(p-1)$-dimensional vector $c_S$ which has
the entries not in $S$ set to zero, and we then let $c_{-S}= c- c_S$.

For a positive semi-definite matrix $A$ we let
$\Lambda_{\rm min}^2 (A)$ be its smallest eigenvalue and
$\Lambda_{\rm max}^2 (A)$ be its largest eigenvalue. 
The smallest eigenvalue of $\Sigma$ is written shorthand as $\Lambda_{\rm min}^2:= \Lambda_{\rm min}^2 (\Sigma)$.
Recall we assume throughout that 
$\Lambda_{\rm min}^2$ stays away from zero: $1/ \Lambda_{\rm min}^2 = {\mathcal O} (1)$.
We use the shorthand notation $\gg 0$ for ``strictly positive and staying away from zero". Thus
throughout we assume $\Lambda_{\rm min}^2 \gg 0 $. 

In order to be able to construct confidence intervals one needs some uniformity in 
unknown parameters. We introduce the following definition. 

\begin{definition} Let ${\cal B}$ be the model for $\beta^0$. 
Let $\{ {\bf Z}_n \}$ be a sequence of real-valued random variables depending on $(X,Y)$
and $\{ r_n \}$ a sequence of positive numbers.
We say that $\{ {\bf Z}_n \}$ is ${\mathcal O}_{\PP_{\beta^0} } ( r_n) $ uniformly in $\beta^0 \in {\cal B}$ if
$$ \lim_{M \rightarrow \infty} \limsup_{n \rightarrow \infty} \sup_{\beta^0 \in {\cal B}}
\PP_{\beta^0} ( | {\bf Z}_n|  > M r_n  ) = 0 . $$
We say that ${\bf Z}_n= o_{\PP_{\beta^0} }(r_n)$ uniformly in $\beta^0 \in {\cal B}$ if 
  $$ \lim_{n \rightarrow \infty} \sup_{\beta_0 \in {\cal B} } \PP_{\beta^0} \biggl (|{\bf Z}_n| > \eta r_n \biggr ) = 0 , \ \forall \ \eta  >0 .$$
  \end{definition}

\subsection{Organization of the rest of the paper}\label{organization.section}
%
Section \ref{knownSigma.section} contains the results for the case $\Sigma$
known and applying a debiased Lasso using
sample splitting. Here we also introduce the concept of an
eligible pair $(\gamma^{\sharp} , \lambda^{\sharp})$ in Definition \ref{gammasharp.condition}. 
Section \ref{eligible.section} contains results and constructions for eligible pairs.
Section \ref{unknownSigma.section} considers the case $\Sigma$ unknown and a debiased Lasso
(without sample splitting). Section \ref{conclusion.section} concludes and Section \ref{proofs.section}
collects the proofs. Section \ref{inequalities.section} (included for completeness)
contains some elementary probability inequalities for products of
Gaussians, which are applied in Section \ref{unknownSigma.section}.

In Tables \ref{conditionsknown.table} and \ref{conditionsunknown.table} we summarize the (sparsity) conditions  we use, see Examples
\ref{l0.example}, \ref{l1.example} and \ref{lr.example} for the case $\Sigma$ known and
 Examples
\ref{l0-unknown.example} and \ref{lr-unknown.example} for the case $\Sigma$ unknown. The
particular cases $r=0$ and $r=1$ follow from the general case $0 \le r \le 1$ when $\Sigma$ is known.
When $\Sigma$ is unknown the case $r=0$ also follows from the general case $0\le r < 1$. With $r=1$ the model
is then too large. We have displayed the extreme cases separately so that the conditions
for these can be read off directly . In particular for $r=0$ one sees
the standard sparsity conditions known from the literature. For $r=1$ ($\Sigma$ known) one sees that
unlike in the other cases there is no logarithmic gap between conditions for asymptotic normality and
asymptotic efficiency. 

 \begin{table}[!htbp]
\centering
\begin{tabular}{c|clclclcl}
 \ & $\Sigma$ & \ \  \  $\Sigma$ &  $\Sigma$    \\
 \ & known & known & known  \\
 \ & ${\cal B} = \{ \| b \|_0^0 \le s \} $ & ${\cal B} = \{ \| b \|_1 \le \sqrt s \} $ & ${\cal B} = \{ \| b \|_r^r \le \sqrt {s^{2-r}} \} $ \\
\hline
 asymp-& $s= o( {n \over  \log p}  ) $ & $s= o( {n \over \log p} ) $  & $s = o ({ n \over \log p} )$\\
 totic & $\lambda^{\sharp} s { \log^{1 \over 2} p } = o(1) $ &$\lambda^{\sharp} \sqrt {ns} = o(1)$& $ \lambda^{\sharp}  { n^{r \over 2} s^{2-r\over 2}\log^{1-r \over 2}  p}=o(1) $  \\
 norma-  & & & \\
 lity & $\lambda^{\sharp} \| \gamma^{\sharp} \|_1 = o(1) $ & $\lambda^{\sharp} \| \gamma^{\sharp} \|_1 = o(1) $ & $  
 \lambda^{\sharp} \| \gamma^{\sharp} \|_1 =o(1)$ \\
\hline
asymp-  & \ & \\
totic & yes & \ \ \ \ yes & yes \\
linearity & & & \\
\hline
asymp-  & \ & \\
totic & $ \| \gamma^{\sharp} \|_0^0 = {\mathcal O} (s) $  & $\| \gamma^{\sharp} \|_1 = {\mathcal O} (\sqrt {ns} ) $ & $\| \gamma^{\sharp} \|_r^r = {\mathcal O} ( {n^{r \over 2} s^{2-r\over 2 } } )  $ \\
efficiency & & & \\
 \end{tabular}
  \caption{The conditions used to prove asymptotic normality, linearity and efficiency when $\Sigma$ is known. 
  Throughout, $(\gamma^{\sharp}, \lambda^{\sharp})$ is required to be an eligible pair
  (see Definition \ref{gammasharp.condition}), i.e. 
    $ \| \Sigma_{-1,-1} (\gamma^{\sharp} - \gamma^0 )  \|_{\infty} \le \lambda^{\sharp} $
    (and $\lambda^{\sharp} \| \gamma^{\sharp} \|_1 \rightarrow 0 $).
    Asymptotic efficiency is established when $\beta^0$ stays away from the boundary of ${\cal B}$.
    In the case ${\cal B} = \{ \| b \|_0^0 \le s \} $ the conditions on $\gamma^{\sharp}$
    for asymptotic efficiency depend on $\beta^0$. }
  \label{conditionsknown.table}
\end{table}

\begin{table}[!htbp]
\centering
\begin{tabular}{c|clclclcl}
 \ & $\Sigma$ & \ \  \  $\Sigma$     \\
 \ & unknown & unknown   \\
 \ & ${\cal B} = \{ \| b \|_0 \le s \} $ & ${\cal B} = \{ \| b \|_r^r \le \sqrt {s^{2-r}} \} $ \\
 \ & & $0 \le r < 1$  \\ 
\hline
 asymp-& $s= o( {\sqrt n \over  \log p}  ) $ & $s= o( {n^{1-r \over 2- r}  / \log p} ) $  \\
 totic & $\lambda^{\sharp} = {\mathcal O} (\sqrt {\log p \over n} )$& $ \lambda^{\sharp}= {\mathcal O} (\sqrt {\log p \over n} )$  \\
 norma-  & & & \\
 lity &  $ \sqrt {\log p\over n} \| \gamma^{\sharp} \|_1 =o(1)$ & $ \sqrt {\log p\over n} \| \gamma^{\sharp} \|_1 =o(1)$ \\
\hline
asymp-  & \ & \\
totic & $\| \gamma^{\sharp}\|_{\rm r}^{\rm r}  = o( {n^{1-{\rm r} \over 2} / \log^{2-{\rm r}\over 2} p} ) $ & $\| \gamma^{\sharp}\|_{\rm r}^{\rm r}  = o( {n^{1-{\rm r} \over 2} / \log^{{2-{\rm r}\over 2} }p} ) $ \\
linearity &  \\
\hline
asymp-  & \ & \\
totic & $ \| \gamma^{\sharp} \|_0^0 = {\mathcal O} (s) $  & $\| \gamma^{\sharp} \|_r^r = {\mathcal O} ( {n^{r\over 2} s^{2-r\over 2 }} ) $ \\
efficiency & & & \\
 \end{tabular}
  \caption{The conditions used to prove asymptotic normality, linearity and efficiency
  when $\Sigma$ is unknown.
  Throughout, $(\gamma^{\sharp}, \lambda^{\sharp})$ is required to be an eligible pair
  (see Definition \ref{gammasharp.condition}), i.e. 
    $ \| \Sigma_{-1,-1} (\gamma^{\sharp} - \gamma^0 )  \|_{\infty} \le \lambda^{\sharp} $
    (and $\lambda^{\sharp} \| \gamma^{\sharp} \|_1 \rightarrow 0 $). The value of ${\rm r}$ may be
    different from $r$. It is assumed to
    be fixed and $0 \le {\rm r} \le 1$. 
    Asymptotic efficiency is established when $\beta^0$ stays away from the boundary of ${\cal B}$.
    In the case ${\cal B} = \{ \| b \|_0^0 \le s \} $ the conditions on $\gamma^{\sharp}$
    for asymptotic efficiency depend on $\beta^0$.}
  \label{conditionsunknown.table}
\end{table}

\vfill\eject
\section{The case $\Sigma$ known}\label{knownSigma.section}

Before presenting ``eligible pairs" in Definition \ref{gammasharp.condition}, we provide the motivation that
led us to this concept.

Recall the  debiased Lasso given in (\ref{de-biased.equation}).
If $\Sigma$ is known we choose the direction $\tilde \Theta_1=: \Theta_1^{\sharp}$ non-random,
 depending on $\Sigma$. 
We invoke the decomposition
$$\hat b_1 - \beta_1^0 =\Theta_1^{\sharp T} X^T \epsilon / n +  \underbrace{({\rm e}_1 -
 \hat \Sigma \Theta^{\sharp} )^T ( \hat \beta - \beta^0 )}_{\rm remainder} . $$
The ${\rm remainder}$  is
$$ ({\rm e}_1- \hat \Sigma \Theta^{\sharp})^T ( \hat \beta - \beta^0 )=
 \underbrace{ \Theta_1^{\sharp T} ( \Sigma - \hat \Sigma)(\hat \beta - \beta^0 )}_{:=(i)} + 
\underbrace{( {\rm e}_1  -  \Sigma \Theta^{\sharp})^T (\hat \beta - \beta^0 )}_{:=(ii)} . 
 $$
 
 The first term $(i)$ can be handled assuming $\Theta_1^{\sharp T} \Sigma \Theta_1^{\sharp} = {\mathcal O} (1)$ and
 $\| \Sigma^{1/2} (\hat \beta - \beta^0)\|_2 = o_{\PP} (1)$. This goes along the lines of
 techniques as in \cite{javanmard2018debias}, applying the conditions
 used there. One then arrives at $(i)= o_{\PP} (1/\sqrt n )$.
 (We will however alternatively use a sample splitting technique later on in 
 Theorem \ref{bsharp.theorem} to simplify the derivations.)
  
 The second term $(ii)$ is additional bias and will be our major concern. If $\Theta_1^{\sharp}= \Theta_1$ this
 term vanishes. However, as we will see it is useful to apply  instead of $\Theta_1$
 some sparse approximation of $\Theta_1$. In fact, we
 aim at a sparse approximation $\Theta_1^{\sharp}$ with  $\Theta_{1,1}^{\sharp}$ being smaller than $\Theta_{1,1}$ and their
 difference not vanishing.  
 
 We will assume conditions
 that ensure the additional bias is negligible and
invoke that
 \begin{equation} \label{dualnorm.equation}
  |( {\rm e}_1  -  \Sigma \Theta_1^{\sharp})^T (\hat \beta - \beta^0 ) | \le 
  \| \Sigma ( \Theta_1^{\sharp} - \Theta_1) \|_{\infty} \| \hat \beta - \beta^0 \|_{1}  
  \end{equation}
   (recall that by the definition of $\Theta_1$ it is true that ${\rm e}_1 = \Sigma \Theta_1$). 
   
   \begin{remark}
 One may think of applying instead  the Cauchy-Schwarz inequality
 $$ |( {\rm e}_1  -  \Sigma \Theta_1^{\sharp})^T (\hat \beta - \beta^0 ) |\le 
 \|  \Sigma^{1/2} ( \Theta_1^{\sharp}- \Theta_1) \|_{2} \| \Sigma^{1/2} (\hat \beta - \beta^0) \|_2  .$$
 This
 leads to requiring that  $\|  \Sigma^{1/2} ( \Theta_1^{\sharp}- \Theta_1) \|_{2} \rightarrow 0$ fast enough. But we actually want
 $\|  \Sigma^{1/2} ( \Theta_1^{\sharp}- \Theta_1) \|_{2} \not\rightarrow 0$ in order to be able to arrive
 at an improvement over the asymptotic variance. 
 
 \end{remark}
 
 \begin{remark}Consider now for some ${\rm p}\ge 1$ the general
 dual norm inequality
 $$ |( {\rm e}_1  -  \Sigma \Theta^{\sharp})^T (\hat \beta - \beta^0 ) |\le 
\| \Sigma ( \Theta_1^{\sharp} - \Theta_1) \|_{\rm p} \| \hat \beta - \beta^0 \|_{\rm q} $$
where $1/ {\rm p} + 1/ {\rm q} =1$. Choosing ${\rm p} \le 2$ here again works against our aim
to improve the asymptotic variance. Thus we need to choose $p>2$ (and therefore
${\rm q} < 2$). This certifies the choice ${\rm p}=\infty$ as being quite natural.
 \end{remark}

  \begin{remark} \label{high-dim.remark} We note that
 $\| \Sigma^{1/2} ( \Theta_1^{\sharp} - \Theta_1) \|_2 \le \| \Sigma (\Theta_1^{\sharp} - \Theta_1) \|_2 / 
 \Lambda_{\rm min}$. This means we want $\| \Sigma (\Theta_1^{\sharp} - \Theta_1) \|_2 \not\rightarrow 0 $ as we
 assume that $\Lambda_{\rm min}$ stays away from zero. Now it is clear that if for some vector
 $v \in \R^p$, it holds that $\| v \|_{\infty} \le \lambda_0^{\sharp}$, then
 $\| v \|_2 \le \sqrt {p} \lambda_0^{\sharp} $. So $\| v \|_2 \not\rightarrow 0 $ implies
 $p\lambda^{\sharp 2}  \not\rightarrow 0  $. In other words, we can only improve the asymptotic variance
 in the high-dimensional case. 
 \end{remark}

Taking the dual norm inequality (\ref{dualnorm.equation}) as starting point we now need
   \begin{equation}\label{boundsharp.equation}
 \| \Sigma (\Theta_1^{\sharp }  - \Theta_1)\|_{\infty} \le \lambda_0^{\sharp} 
 \end{equation}
 for some constant $\lambda_0^{\sharp}$ small enough, such that uniformly in $\beta^0 \in {\cal B}$
 $$\lambda_0^{\sharp} \| \hat \beta - \beta^0 \|_1 = o_{\PP_{\beta^0}} (1/\sqrt n  ) .$$

With the above considerations as motivation, we
 now concentrate on 
 an $\ell_{\infty}$-condition as given in  inequality (\ref{boundsharp.equation}). 
 We fix some $\lambda^{\sharp}$ and construct vectors $\Theta_1^{\sharp}$
 for which inequality (\ref{boundsharp.equation}) holds. 
 It is based on replacing the vector of coefficients $\gamma^0$ of the regression
 of ${\bf x}_{1}$ on ${\bf x}_{-1}$ by a sparse approximation $\gamma^{\sharp}$.
 
 \begin{definition}\label{gammasharp.condition}
 Let $\gamma^{\sharp} \in \R^{p-1}$ and $\lambda^{\sharp} >0$.
 We say that the pair $(\gamma^{\sharp} , \lambda^{\sharp} )$ is eligible if 
 \begin{equation}\label{linfty.equation}
 \|  \Sigma_{-1,-1} (\gamma^{\sharp}- \gamma^0) \|_{\infty} \le \lambda^{\sharp} 
 \end{equation}
 and
 \begin{equation}\label{l1sparse.equation}
  \lambda^{\sharp} \| \gamma^{\sharp} \|_1 \rightarrow 0 .
  \end{equation}
   \end{definition}
 
 The two conditions in Definition \ref{gammasharp.condition} will allow us to arrive at
 (\ref{boundsharp.equation}) as is shown in the next lemma.
  
 \begin{lemma} \label{need.lemma} 
 Suppose $(\gamma^{\sharp}, \lambda^{\sharp})$ is an eligible pair.
    Then
  $$1- \gamma^{0T} \Sigma_{-1,-1} \gamma^{\sharp} \ge \Lambda_{\rm min}^2 -o(1)$$
  i.e., $1- \gamma^{0T} \Sigma_{-1,-1} \gamma^{\sharp}  \gg 0 $ eventually. 
  Let (for $n$ sufficiently large)
  \begin{equation} \label{Thetasharp.equation}
 \Theta_1^{\sharp} := \begin{pmatrix} 1 \cr - \gamma^{\sharp} \end{pmatrix}/ (1- \gamma^{0T} \Sigma_{-1, -1} \gamma^{\sharp} ).
 \end{equation} 
 Then we have
 $$\| \Sigma( \Theta_1^{\sharp} - \Theta_1)  \|_{\infty} \le \lambda_0^{\sharp}$$
 where 
 $$\lambda_0^{\sharp}:= \lambda^{\sharp} /  (1- \gamma^{0T} \Sigma_{-1, -1} \gamma^{\sharp}) =
 {\mathcal O} (\lambda^{\sharp} ). $$
 Moreover,
  \begin{eqnarray*}
    \Theta_1^{\sharp T } \Sigma \Theta_1^{\sharp}&=& \Theta_{1,1}^{\sharp} +o(1)\\
 & \le & \Theta_{1,1} + o(1).
 \end{eqnarray*}
Finally, in order to have a non-vanishing  improvement  of
$\Theta_{1,1}^{\sharp}$ over $\Theta_{1,1}$ it must be true that
 $\gamma^0$ is not sparse, in the sense that
 $$ \lambda^{\sharp} \| \gamma^0 \|_1\gg 0 .$$
  \end{lemma} 
  
  \begin{remark}\label{regression.remark}
  The first condition (\ref{linfty.equation}) of Definition 
  \ref{gammasharp.condition} can be rewritten as
 $$ {\bf x}_{-1} \gamma^0 = {\bf x}_{-1} \gamma^{\sharp} + \xi^0, \ | \EE {\bf x}_j \xi^0 | \le
 \lambda^{\sharp}  \ \forall \ j \in \{2 , \ldots , p \} .$$
 The second condition (\ref{l1sparse.equation}) in this definition can be thought of as a sparsity condition on $\gamma^{\sharp} $. 
 The two conditions together require that the regression of
${\bf x}_1$ on ${\bf x}_{-1}$ is sparse when one relaxes the orthogonality condition
of residuals to approximate orthogonality.\footnote{Condition (\ref{l1sparse.equation}) is of the same nature
as the condition $\lambda \| \beta^0 \|_1 = o(1)$ (which follows from
the classical condition $\lambda^2 \| \beta^0 \|_0^0 =o(1)$ if
$\| \beta^0 \|_2 = {\mathcal O} (1)$) when applying the Lasso (\ref{Lasso.equation})
with tuning parameter $\lambda$.}  One may think of $\gamma^0$ as a ``least squares estimate"
of $\gamma^{\sharp}$ in a noisy regression model. We refer to Subsection \ref{regression.section}
which gives a very natural interpretation of eligible pairs. 
 Further, for $\Theta_1^{\sharp}$ defined in (\ref{Thetasharp.equation})
 one has the equivalence
 $$ \Theta_{1,1} - \Theta_{1,1}^{\sharp}\gg 0 \ \Leftrightarrow \ \EE ( \xi^0)^2 \gg 0 . $$
 \end{remark}

  To have $\Theta_{1,1}^{\sharp}$ improving over $\Theta_{1,1}$ we see from the above lemma that we aim at a situation where $\gamma^0$, and hence $\Theta_1$,  is not sparse,
  but where $\gamma^0$ can be replaced by a sparse vector $\gamma^{\sharp}$. 
  For some special $\Sigma$'s, we give examples of eligible pairs in Section \ref{eligible.section}.
  That section also discusses for a given $\lambda^{\sharp}$ uniqueness of the vectors $\gamma^{\sharp}$
  for which the pair $(\gamma^{\sharp}, \lambda^{\sharp})$ is eligible. Moreover, we show cases where
  ${\bf x}_{-1} \gamma^{\sharp}$ is an approximation of the projection of ${\bf x}_1$ on a subset of
  ${\bf x}_S $ of the other variables for some $S \subset \{ 2 , \ldots , p \} $, see Lemma \ref{projectS.lemma}.
  This may help to understand why the Cram\'er Rao lower bound can be achieved in those
  cases. 
      
  Lemma \ref{need.lemma} essentially has all the ingredients to prove asymptotic normality
  of the debiased Lasso (\ref{de-biased.equation}) with direction $\tilde \Theta_1= \Theta_1^{\sharp}$ and
  $\Theta_1^{\sharp}$ given in (\ref{Thetasharp.equation}) in this lemma. 
It can be done along the lines of Theorem 3.8 in \cite{javanmard2018debias}, assuming the conditions
stated there. However, as the authors point out, 
when using instead the sample splitting approach their Assumption $(iii)$ is not needed. 
Sampling splitting it is also mathematically less involved. We present it here.

Assume the sample size $n$ is even.
Define the matrices
\begin{eqnarray*}
(X_I, Y_I) &:= &\{X_{i,1}, \ldots , X_{i,p}, Y_i  \}_{1 \le i \le n/2} \in \R^{n/2 \times (p+1)}, \\
 (X_{II}, Y_I) &:=& \{X_{i,1}, \ldots , X_{i,p}, Y_i  \}_{n/2 < i \le n}  \in \R^{n/2 \times (p+1)}. 
\end{eqnarray*}
Let $\hat \beta_I$ be an estimator of $\beta^0$ based on the
first half $(X_I, Y_I)$ of the sample, for instance the Lasso estimator
$ \arg \min  \{ \| Y_I - X_I b\|_2^2 /n+  \lambda \| b \|_1 \}$. 
Similarly, let $\hat \beta_{II}$ be an estimator of $\beta^0$ based on the second half $(X_{II}, Y_{II})$ of the sample.
Let $(\gamma^{\sharp}, \lambda^{\sharp})$ be an eligible pair.
We then define the two debiased estimators
\begin{eqnarray*}\label{two-sample.equation}
\hat b_{I,1}^{\sharp} &:=& \hat \beta_{II,1} + 2\Theta_1^{\sharp T} X_I^T\biggl (Y_I - X_I \hat \beta_{II} \biggr )/n\\
\hat b_{II,1}^{\sharp} &:=& \hat \beta_{I,1} + 2\Theta_1^{\sharp T} X_{II}^T\biggl (Y_{II} - X_{II}\hat \beta_I \biggr )/n
\end{eqnarray*}
where $\Theta_1^{\sharp} $ is given in (\ref{Thetasharp.equation}) in Lemma \ref{need.lemma}.
The final estimator $\hat b_1^{\sharp}$ is obtained  by averaging these two:
\begin{equation}\label{bsharp.equation}
\hat b_1^{\sharp} := {\hat b_{I,1}^{\sharp}+\hat b_{II,1}^{\sharp} \over 2} . 
\end{equation}

Let now ${\cal B}$ be a given model class for the unknown vector of regression
coefficients $\beta^0$.

\begin{theorem} \label{bsharp.theorem}
Let $(\gamma^{\sharp}, \lambda^{\sharp})$ be an eligible pair,
$\Theta_1^{\sharp} $ be given in (\ref{Thetasharp.equation}) and 
$\hat b_1^{\sharp}$ be given in 
(\ref{bsharp.equation}) with $\hat b_{I,1}^{\sharp}$ and $\hat b_{II,1}^{\sharp}$ the debiased estimators
based on $\Theta_1^{\sharp}$ using the splitted sample.
Suppose that uniformly in $\beta^0 \in {\cal B}$
\begin{equation}\label{l2-convergence.equation}
 \| \Sigma^{1/2} (\hat \beta_I - \beta^0) \|_2 = o_{\PP_{\beta^0}} (1) ,\ \| \Sigma^{1/2} (\hat \beta_{II} - \beta^0 )\|_2 = o_{\PP_{\beta^0}} (1) 
 \end{equation}
and
\begin{equation}\label{l1-convergence.equation}
\sqrt n \lambda^{\sharp} \| \hat \beta_I - \beta^0 \|_1 = o_{\PP_{\beta^0}} (1) ,\ 
\sqrt n \lambda^{\sharp} \| \hat \beta_{II} - \beta^0 \|_1 = o_{\PP_{\beta^0}} (1) . 
\end{equation}
Then, uniformly in $\beta^0 \in {\cal B}$,
$$ \hat b_1^{\sharp} - \beta_1^0 = \Theta_1^{\sharp T } X^T \epsilon /n + o_{\PP_{\beta^0} } (1/ \sqrt n ) , $$
and 
$$ \lim_{n \rightarrow \infty} \sup_{\beta_0 \in {\cal B} } \PP_{\beta^0}\biggl (
{\sqrt n ( \hat b_1^{\sharp} - \beta_1^0 ) \over \Theta_{1,1}^{\sharp} } \le z \biggr )= \Phi(z) , \ \forall \ z \in \R. $$
\end{theorem} 

\begin{corollary} \label{efficient.corollary} 
Theorem \ref{bsharp.theorem} shows that under its conditions
the estimator $\hat b_1^{\sharp}$
is uniformly asymptotically linear and regular.
It means that for this estimator the Cram\'er Rao lower bound as
given in Subsection \ref{CRLB.section} is relevant.
Depending among other things on ${\cal B}$ and $\Sigma$ it does or does not reach the Cram\'er Rao lower bound,
see Remark \ref{lemmas.remark}.
\end{corollary}

\begin{remark} \label{lemmas.remark} Lemmas \ref{direct.lemma}, \ref{reversed-irrepresentable.lemma} and
\ref{non-sparse.lemma} in Subsection \ref{engineering.section} present examples of eligible pairs
$(\gamma^{\sharp}, \lambda^{\sharp})$ where $\Theta_{1,1} - \Theta_{1,1}^{\sharp}$ 
(with $\Theta_{1,1}^{\sharp}$ given in (\ref{Thetasharp.equation})) is
non-vanishing:  $\Theta_{1,1} - \Theta_{1,1}^{\sharp} \gg 0 $.Thus for those cases
Theorem \ref{bsharp.theorem} shows an asymptotic variance remaining strictly
smaller than $\Theta_{1,1} $. Moreover, the constructions 
of Lemmas \ref{direct.lemma}, \ref{reversed-irrepresentable.lemma} and
\ref{non-sparse.lemma} allow for directions $\Theta_1^{\sharp} $ 
(depending on $ \| \beta^0 \|_0^0\in  {\cal H}_{\beta^0}$
in model (\ref{sparse-model.equation}). In models
(\ref{l1-model.equation}) and (\ref{lr-model.equation}) the direction
$\Theta_1^{\sharp}$ lies in ${\cal H}_{\beta^0 }$ after scaling where only
the scaling depends on $\beta^0$.
For these cases (special constructions of $\Sigma$) the Cram\'er Rao lower bound is therefore achieved.
\end{remark} 

We now discuss the requirements (\ref{l2-convergence.equation}) and
(\ref{l1-convergence.equation}) for the models (\ref{sparse-model.equation}),
(\ref{l1-model.equation}) and (\ref{lr-model.equation}). The overall picture
is summarized in Table \ref{conditionsknown.table}. 

\begin{example}\label{l0.example}
Consider model 
model (\ref{sparse-model.equation}) with $s= o( n / \log p )$ (sparsity variant (i)).
For the 
Lasso estimator 
$\hat \beta $ given in (\ref{Lasso.equation}),
with appropriate choice of the tuning parameter $\lambda \asymp \sqrt {\log p / n } $, one has
uniformly in $\beta^0 \in {\cal B}$
$$\|\Sigma^{1/2} (\hat \beta- \beta^0 )\|_2^2 = {\mathcal O}_{\PP_{\beta^0}} ( s {  \log p / n}  ) , \
\| \hat \beta - \beta^0 \|_1 = O_{\PP_{\beta^0}} ( s \sqrt {  \log p / n}  ). $$
This follows from e.g.\ Theorem 6.1 in \cite{BvdG2011},  together with results for Gaussian
quadratic forms as given in (\ref{quadratic-forms.equation}). 
So for $\hat \beta_I$ and $\hat \beta_{II}$ Lasso's based on the splitted sample
with suitable tuning parameter $\lambda$, the requirement on $\lambda^{\sharp}$ becomes
$$ \lambda^{\sharp} s \sqrt {\log p}= o( 1).  $$
This is guaranteed  when $ \lambda^{\sharp}= {\mathcal O} ( \sqrt {\log p} /n )$.
The smaller the sparsity $s$, the more room there is for improvement (i.e., the larger
the collection
of covariance matrices $\Sigma$ that allow for improvement over $\Theta_{1,1}$).
For example, if in fact
$s= o( \sqrt {n } / \log p ) $
we can take
$ \lambda^{\sharp}= {\mathcal O} ( \sqrt {\log p /n} )$.
Finally, if $\| \gamma^{\sharp} \|_0^0 = s_{\sharp}$ with $s_{\sharp} \le s- s_0 -1$ the estimator
$\hat b_1^{\sharp}$ is asymptotically efficient. 
\end{example}

\begin{remark} \label{efficient.remark} 
In view of Theorem \ref{bsharp.theorem} and the statements of Example \ref{l0.example} and Remark
\ref{lemmas.remark}, we see that for model
(\ref{sparse-model.equation}) the debiased
estimator of Theorem 3.8 in \cite{javanmard2018debias}, which uses
$\tilde \Theta_1= \Theta_1$, is in certain cases asymptotically  inefficient as
a choice $\tilde \Theta_1 = \Theta_1^{\sharp}\not= \Theta_1$ can give
an improvement in the asymptotic variance (and is then efficient for certain
cases).
\end{remark}

\begin{example}\label{l1.example}
In this example we take the model (\ref{l1-model.equation}) with $0< s = o( n / \log p ) $. 
Let $\hat \beta$ be again the Lasso estimator 
given in (\ref{Lasso.equation})
with appropriate choice of the tuning parameter $\lambda \asymp \sqrt {\log p / n } $.
One may use a ``slow rates" result: uniformly in $\beta^0 \in {\cal B}$ it is true that
$$ \| \Sigma^{1/2} ( \hat \beta - \beta^0 ) \|_2 = o_{\PP_{\beta^0} } (1) , \ 
\| \hat \beta - \beta^0 \|_1 = {\mathcal O}_{\PP_{\beta^0 } } (\sqrt s) , $$
see for example
\cite{BvdG2011}, Theorem 6.3, and combine this with results for quadratic forms as
given in (\ref{quadratic-forms.equation}). 
(The arguments for establishing these ``slow rates" are in fact as in Lemma \ref{slowLasso.lemma} for the node-wise Lasso.) Taking $\hat \beta_I$ and $\hat \beta_{II}$ again as appropriate Lasso's, condition (\ref{l1-convergence.equation}) 
of Theorem \ref{bsharp.theorem} holds if
$$\lambda^{\sharp} \sqrt {ns}= 
o(1) . $$
Then for $\| \gamma^{\sharp} \|_1 = {\mathcal O} (\sqrt {ns})$ we get an eligible
pair $( \gamma^{\sharp}, \lambda^{\sharp})$. Since when $\beta^0$ stays away from the 
boundary, such $\gamma^{\sharp}$ form a model direction
after scaling, we see that the Cram\'er Rao lower bound is achieved.
Recall now that in Example \ref{l1-model.example} we concluded that, in view
of Corollary \ref{CRLB.corollary},
in order to be able to improve over $\Theta_{1,1}$ we must have $\| \gamma^0 \|_1 $
of order larger than $\sqrt {ns} $. 
As pointed out in Remark \ref{lemmas.remark} 
one may apply one of
the  Lemmas \ref{direct.lemma}, \ref{reversed-irrepresentable.lemma} or
\ref{non-sparse.lemma} to create vectors $\gamma^0$ and eligible pairs
$(\gamma^{\sharp}, \lambda^{\sharp})$ where $\lambda^{\sharp} \sqrt {ns} 
\rightarrow 0$, $\| \gamma^{\sharp} \|_1$ can take any value of order $ \sqrt {ns} $ and 
$\Theta_{1,1} - \Theta_{1,1}^{\sharp} \gg 0$.
\end{example}

\begin{example} \label{lr.example} 
The results for model (\ref{l1-model.equation})
are a special case of those for model (\ref{lr-model.equation}).
One needs again $s= o( n / \log p)$ and one may for example apply
Corollary 2.4 in \cite{vdG2016}, again together with (\ref{quadratic-forms.equation}). For the Lasso
$\hat \beta$ in (\ref{Lasso.equation}) with $\lambda \asymp \sqrt {\log p / n}$ appropriately chosen one gets
uniformly in $\beta^0 \in {\cal B}$
$$ \| \Sigma^{1/2} ( \hat \beta - \beta^0 ) \|_2^2 = {\mathcal O}_{\PP_{\beta^0}} 
( {s \log p / n  } )^{2-r \over 2 }, \ \| \hat \beta - \beta^0 \|_1 =
{\mathcal O}_{\PP_{\beta^0 }} ( ( \log p / n)^{1-r \over 2 } s^{2-r \over 2} ). $$
The requirement (\ref{l1-convergence.equation}) on $\lambda^{\sharp}$ is thus
$$ \lambda^{\sharp} ( \log p)^{1- r \over 2} \sqrt {n^{r } s^{2-r } }= o(1) . $$
If $\| \gamma^{\sharp} \|_r^r = {\mathcal O} ( \sqrt { n^r s^{2-r} } )$
then $(\gamma^{\sharp}, \lambda^{\sharp})$ is an eligible pair and the Cram\'er Rao
lower bound is achieved whenever $\beta^0$ stays away from the boundary. 
In order to be able to improve over $\Theta_{1,1}$ we now
need $\| \gamma^0 \|_r^r $ of larger order $\sqrt { n^r s^{2-r} }$
by Corollary \ref{CRLB.corollary}. 
Remark \ref{lemmas.remark} can be taken into the considerations here too.

\end{example}

\section{Finding eligible pairs} \label{eligible.section}
The main results of this section can be found in Subsection
\ref{engineering.section} where for any $\lambda^{\sharp}$ we construct eligible pairs
$(\gamma^{\sharp}, \lambda^{\sharp})$ by choosing $\gamma^0$ appropriately.
The results can be seen as existence proofs. Before doing these
constructions we discuss uniqueness in Subsection
\ref{uniqueness.section}, in Subsection \ref{usingLasso.section} the noiseless Lasso as a practical method
for improving over $\Theta_{1,1}$ and in Subsection \ref{usingprojections.section} we examine
whether or not projections on a subset
of the variables can lead to eligible pairs. For the latter we impose rather string conditions.
We show in Subsection \ref{approxprojections.section} that eligible pairs are more flexible than projections.
Nevertheless in the final part of this section we return to projections as they come
up naturally when imposing non-sparsity constraints on $\gamma^0$

\subsection{Uniqueness}\label{uniqueness.section}

\begin{lemma}\label{unique.lemma} Fix $\lambda^{\sharp}$ and let
$(\gamma^{\sharp}, \lambda^{\sharp})$ and $(\gamma^{\flat}, \lambda^{\sharp})$
be two eligible pairs. Then
$$ \|  \Sigma_{-1,-1}^{1/2}  (\gamma^{\sharp} - \gamma^{\flat} )\|_2^2 \rightarrow 0 . $$
\end{lemma}

This lemma tells us that asymptotically it makes no difference to debias using
$(\gamma^{\sharp}, \lambda^{\sharp})$ or $(\gamma^{\flat}, \lambda^{\sharp})$. 


\subsection{Using the Lasso}\label{usingLasso.section}
Fix some tuning parameter $\lambda_{\rm Lasso}$ and consider the noiseless Lasso
\begin{equation}\label{noiselessLasso.equation}
\gamma_{\rm Lasso}:= \arg \min_{c \in \R^{p-1}}
\biggl \{ \EE ( {\bf x}_1 - {\bf x}_{-1} c )^2 + 2  \lambda_{\rm Lasso} \| c \|_1 \biggr \} . 
\end{equation}
One may verify that  $(\gamma_{\rm Lasso}, \lambda_{\rm Lasso})$ is an eligible pair if
$ \lambda_{\rm Lasso} \| \gamma^0 \|_1 \rightarrow 0 $. But the latter is exactly
what we want to avoid. 
If $( \gamma^{\sharp} , \lambda^{\sharp})$ is an eligible pair
then the noiseless Lasso can find it if one chooses $\lambda_{\rm Lasso}= {\mathcal O} ( \lambda^{\sharp})$ larger
than $\lambda^{\sharp}$, as follows from the next lemma. It says that given $\lambda^{\sharp}$ one
may use the noiseless Lasso for constructing a direction, $\Theta_{1, {\rm Lasso}}$ say,
with which one has same improvement over $\Theta_{1,1}$ as with $\Theta_1^{\sharp}$.

\begin{lemma}\label{findpair.lemma} Suppose $( \gamma^{\sharp}, \lambda^{\sharp})$ is an eligible pair.
Let $\lambda_{\rm Lasso} > \lambda^{\sharp}$ 
and $ \lambda_{\rm Lasso} \| \gamma^{\sharp} \|_1$ $ \rightarrow 0$.
Let $\gamma_{\rm Lasso}$ be the noiseless Lasso defined in (\ref{noiselessLasso.equation}).
Then
$$\|  \Sigma_{-1,-1}^{1/2}  (\gamma_{\rm Lasso} - \gamma^{\sharp} )\|_2^2 \rightarrow 0 
 . $$
In addition, if $\lambda_{\rm Lasso} \ge 2 \lambda^{\sharp} $ (say)
then
$( \gamma_{\rm Lasso} , \lambda_{\rm Lasso})$ is an eligible pair.
\end{lemma}

\subsection{Using projections}\label{usingprojections.section}

In this subsection we investigate (rather straightforwardly) conditions such that the coefficients of
a projection of ${\bf x}_{-1} \gamma^0$ 
can be joined with a $\lambda^{\sharp}$ to form an eligible pair.

Fix some set $S\subset \{ 2, \ldots , p \}$ with cardinality ${\rm s}$. 
The value of ${\rm s}$ need not be $s$, where $s$ is the sparsity
used in the model class ${\cal B}$.
Let ${\bf x}_S \gamma_S^S $ be the projection
of ${\bf x}_{-1} \gamma^0$ on ${\bf x}_S $.

Let $\gamma^S \in \R^{p-1}$ be the vector $\gamma_S^S
\in \R^{\rm s}$ completed with zeroes. 
Then $\| \gamma^S \|_1 \le \sqrt {\rm s} \| \gamma^S \|_2 = {\mathcal O} (\sqrt {\rm s} )$
so that when 
\begin{equation} \label{lambdasqrts.equation}
\lambda^{\sharp} \sqrt {\rm s} \rightarrow 0
\end{equation}
we have for $\gamma^S$ the sparsity condition (\ref{l1sparse.equation})
of Definition \ref{gammasharp.condition}:
$\lambda^{\sharp} \| \gamma^S \|_1 \rightarrow 0 $. 
Note that $\gamma^{0T} \Sigma_{-1,-1} \gamma^S = \gamma^{ST} \Sigma_{-1,-1} \gamma^S(= \EE ( {\bf x}_S \gamma_S^S )^2)$ exactly in this case.
Furthermore
$$ v^S := \Sigma_{-1,-1} ( \gamma^0 - \gamma^S) = \begin{pmatrix}  0 \cr  v_{-S}^S  \end{pmatrix} $$
where
$$ v_{-S}^S = \EE {\bf x}_{-S} {\bf x}_{-1} \gamma^{-S} , $$
with
$$ {\bf x}_{-1} \gamma^{-S} = ({\bf x}_{-S} \gamma_{-S}) {\rm A} {\bf x}_S =:
\biggl ( {\bf x}_{-S} - {\bf x}_S^T \Sigma_{S,S}^{-1} \Sigma_{S,-S} \biggr ) \gamma_{-S}^0 $$
being the anti-projection of ${\bf x}_{-S} \gamma_{-S}^0$ on ${\bf x}_S$. 
We check whether
the pair $(\gamma^S , \lambda^{\sharp})$ is an eligible
pair, which is the case if  $\lambda^{\sharp} \sqrt {\rm s} = o(1)$ and $\| v_{-S}^S \|_{\infty} \le \lambda^{\sharp} $.  
We briefly discuss some conditions that may help ensuring the latter.

Let $\vvert A \vvert_1 := \max_{j } \sum_{k} |a_{j,k} | $
be the $\ell_{1}$-operator norm of the matrix  $A $. 

\begin{lemma} \label{antiprojection.lemma} It holds that
$$ \| v_{-S}^S \|_{\infty} \le   \vvert  \Sigma_{-S, -S} - \Sigma_{-S, S} \Sigma_{S,S}^{-1} \Sigma_{S,-S}  \vvert_1   \| \gamma_{-S}^0 \|_{\infty}  . $$

\end{lemma}

To arrive at a sparse approximation of $\gamma^0$ one may consider putting its
$p-{\rm s}-1$ smallest-in-absolute-value coefficients to zero.
To this end, consider the ordered sequence $| \gamma^0 |_{(1)} \ge \cdots \ge |\gamma^0|_{(p-1)} $, 
write $|\gamma^0|_{(j)} =: |\gamma_{r_j }^0 |$, $j= 1, \ldots , p-1 $, and let
$S= \{ r_1 , \ldots r_{\rm s} \} $. Then clearly $\| \gamma_{-S}^0 \|_{\infty} \le \| \gamma^0 \|_2 / \sqrt {\rm s}
= {\mathcal O} (1/ \sqrt {s}) $.
The condition (\ref{lambdasqrts.equation}) excludes $1/ \sqrt {\rm s} = {\mathcal O} ( \lambda^{\sharp})$.
We therefore examine situations where the coefficients in $\gamma^0$ decrease at a rate quicker than $1/ \sqrt {\rm s}$

\begin{definition} Let $N$ be some integer and $0 \not= v \in \R^N$ be a vector.  We call
$$\kappa (v) := {  \| v \|_1 \| v \|_{\infty} \over \| v \|_2 }  $$
the sparsity index of $v$.
We say that $v$ is asymptotically sparse if $\kappa (v)\rightarrow 0$.
\end{definition}

A vector $v$ with $\| v \|_2={\mathcal O} (1)$ can have some 
relatively large coefficients, but it cannot have too many of these. If 
in addition $\| v \|_1$ is large it cannot have many zeroes either.
Asymptotic non-sparseness of the vector $v$ with the large
coefficients removed means that there are many very small
non-zero coefficients. 

\begin{example} Let $v_j = 1/ \sqrt {N} $ for all $j$. Then 
$\| v \|_1 = 1/ \| v\|_{\infty} = \sqrt N  $. Thus $\kappa (v)=1$ 
and $v$ is not (asymptotically) sparse.
\end{example}

\begin{example} Let ${\rm s} \le N$ and $v_j=0$ for $j \le {\rm s}$ and $v_j = 1/ \sqrt {j \log N } $ for $j>{\rm s}$. Then 
$\kappa(v)  \asymp \sqrt{ N /({\rm s}  \log^2 N) }$.
The vector $v$ is not asymptotically sparse if ${\rm s}= {\mathcal O}  (N/ (\log^2 N) )$. 
\end{example}

\begin{lemma} \label{non-sparse1.lemma}
Suppose that $\gamma_{-S}^0$ is not asymptotically sparse.
Let
$$\lambda^S = {\kappa (\gamma_{-S}^0) \over C  \| \gamma_{-S}^0 \|_1 } , $$ where
$$C:= \vvert \Sigma_{-S, -S} - \Sigma_{-S, S} \Sigma_{S,S}^{-1} \Sigma_{S,-S}\vvert_1  . $$
Assume that $\lambda^S \sqrt {\rm s} \rightarrow 0$. Then
$( \gamma^S , \lambda^S)$ is an eligible pair and
$\lambda^S \| \gamma^0 \|_1 \not\rightarrow 0 $.
\end{lemma}

\subsection{Approximate projections}\label{approxprojections.section}
Recall that the first condition (\ref{linfty.equation}) of Definition \ref{gammasharp.condition} can be written as
  $$ {\bf x}_{-1} \gamma^0 = {\bf x}_{-1} \gamma^{\sharp} + \xi^0 $$
  where $\xi^0$ is ``noise" satisfying $\| \EE {\bf x}_{-1}^T \xi^0 \|_{\infty} \le \lambda^{\sharp}$. 
  Denote the active set of $\gamma^{\sharp}$ 
  by $S = \{ j\in \{ 2 , \ldots , p \}:\ \gamma_j^{\sharp} \not= 0 \} $ and its cardinality by ${\rm s} := |S|$.
   Then $\lambda^{\sharp}\sqrt {\rm s} \rightarrow 0$ implies that $1/ \Theta_{1,1}^{\sharp}$ is asymptotically
  the squared residual of the projection of ${\bf x}_1$ on ${\bf x}_S$ as is shown in the next lemma.
  In that sense, eligible pairs are more flexible than projections. 

   If ${\rm s}$ is small enough, then for model (\ref{sparse-model.equation}),
   (\ref{l1-model.equation}) or (\ref{lr-model.equation}) 
  one has that $\Theta_1^{\sharp}$ is a model direction (after rescaling). The estimator
  $\hat b_1^{\sharp}$ in Theorem \ref{bsharp.theorem} then reaches the Cram\'er Rao lower bound.

  \begin{lemma} \label{projectS.lemma} Suppose for some set $S$ 
  with cardinality ${\rm s}$ that $\gamma^{\sharp} = \gamma_S^{\sharp}$
  and assume
 the first condition (\ref{linfty.equation}) of Definition \ref{gammasharp.condition}. If $\lambda^{\sharp} \sqrt {\rm s} \rightarrow 0 $
 then the second condition (\ref{l1sparse.equation}) holds too so that $(\gamma^{\sharp}, \lambda^{\sharp})$ is an eligible pair. Moreover,
 $$1/ \Theta_{1,1}^{\sharp} = \EE({\bf x}_{1} - {\bf x}_S \gamma_S^S )^2 + o(1) $$
 where 
 $$ \gamma_S^S = \Sigma_{S,S}^{-1}  \EE  {\bf x}_S^T {\bf x}_1 .$$
  \end{lemma}

\subsection{Reverse engineering}\label{engineering.section}
In this subsection we fix $\lambda^{\sharp}$ and then construct vectors $\gamma^0 \in \R^{p-1}$ such that there is an eligible pair
$(\gamma^{\sharp}, \lambda^{\sharp} )$. These constructions are in a sense equivalent
but approach the problem from different angles.
In these constructions the vector $\gamma^{\sharp}$
is having active set $S= \{ \gamma_j^{\sharp} \not= 0 \}$ with cardinality ${\rm s}:= |S|$.
The sparsity of $\gamma^{\sharp}$ is then measured in terms of 
the value of ${\rm s}$. More general constructions are possible, but in this
way we can apply the results to any of the models
(\ref{sparse-model.equation}), (\ref{l1-model.equation}) or (\ref{lr-model.equation}).
In view of Lemma \ref{projectS.lemma} is means that the constructions correspond
to an approximate projection on ${\bf x}_S$. 

Throughout this subsection, the matrix $\Sigma_{-1, -1}$ is
is assumed to have 1's one the diagonal and smallest eigenvalue $\Lambda_{\rm min}^2 (\Sigma_{-1, -1})\gg 0$.

\subsubsection{Which $\gamma^0$'s are allowed?}
We let for $\gamma^0 \in \R^p$
$$ \Sigma (\gamma^0):= \begin{pmatrix} 1 & \gamma^{0T} \Sigma_{-1,-1} \cr
\Sigma_{-1, -1} \gamma^0 & \Sigma_{-1, -1} \cr \end{pmatrix} .$$

\begin{definition} \label{allowedgamma.definition} We say that the vector $\gamma^0 $ is allowed if
$\Sigma (\gamma^0) $ is positive definite, with
$ \Lambda_{\rm min}^2 (\Sigma (\gamma^0)) \gg 0 $
and  $\| \Sigma_{-1, -1} \gamma^0 \|_{\infty} \le 1 $. 
\end{definition}

\begin{lemma} \label{allowedgamma.lemma}
Suppose that 
$$1- \| \Sigma_{-1,-1}^{1/2}  \gamma^0 \|_2^2  \gg 0 . $$
Then
$\gamma^0$ is allowed. 
\end{lemma}

\subsubsection{Regression: $\gamma^0$ as least squares estimate of $\gamma^{\sharp}$}\label{regression.section}
In this subsection we create $\gamma^0$ using random noise. We then arrive at an eligible
pair ``with high probability". Let $N\in \Nat$ be a given sequence with $N > p$. Take a matrix
$Z_{-1} \in \R^{N \times (p-1)}$ which has $Z_{-1}^T Z_{-1} / N = \Sigma_{-1,-1} $. 
Let $\xi \in \R^N$ have i.i.d.\ standard Gaussian entries, let
$\gamma^{\sharp} \in \R^{p-1}$ be a vector satisfying $\| \gamma^{\sharp} \|_1= o ( \sqrt {\log p / N } )$ 
 and define 
 $$Z_1 := Z_{-1} \gamma^{\sharp} + \xi . $$
 Let
 $$ \gamma^0 := (Z_{-1}^T Z_{-1} )^{-1} Z_{-1}^T Z_1 $$
 be the least  squares estimator of $\gamma^{\sharp}$. Finally let $\lambda^{\sharp} \asymp \sqrt {\log p / N}$.
 be appropriately chosen. 
 Then $( \gamma^{\sharp} , \lambda^{\sharp})$ is with high probability an eligible pair. Indeed
 the first condition (\ref{linfty.equation}) of Definition \ref{gammasharp.condition} follows from
 $$\|  \Sigma_{-1,-1} (\gamma^0- \gamma^{\sharp} ) \|_{\infty}= 
 \| Z_{-1}^T \xi/N \|_{\infty} = {\mathcal O}_{\PP} ( \sqrt {\log p / N} )$$
 so that for appropriate $\lambda^{\sharp} \asymp \sqrt {\log p / N}$, with high probability
  $$\|  \Sigma_{-1,-1} (\gamma^0- \gamma^{\sharp} ) \|_{\infty}\le \lambda^{\sharp} . $$
 The second condition (\ref{l1sparse.equation}) of Definition \ref{gammasharp.condition} follows from the condition
 $\| \gamma^{\sharp} \|_1= o ( \sqrt {\log p / N } )$ so $\lambda^{\sharp} \| \gamma^{\sharp} \|_1 = o(1)$. 
 We also see that if $p/N \gg 0$ then with high probability $\Theta_{1,1}^{\sharp}$ as given
 in (\ref{Thetasharp.equation}) is an improvement over $\Theta_{1,1}$, since
 $$ \| \Sigma_{-1,-1 }^{1/2} ( \gamma^0 - \gamma^{\sharp} )\|_2^2 = \chi_{p-1}^2 / N $$
 where $\chi_{p-1}^2 $ has the chi-squared distribution with $p-1$ degrees of freedom. 
 It follows that
$$  \| \Sigma_{-1,-1 }^{1/2} ( \gamma^0 - \gamma^{\sharp} )\|_2^2= {p + {\mathcal O} (\sqrt p ) \over N} . $$
With high probability this stays away from zero so that also
$$ \Theta_{1,1} - \Theta_{1,1}^{\sharp} \gg 0 . $$
For appropriate $N$ with $1- p/N \gg 0$ the vector $\gamma^0$ is with high probability also allowed by
Lemma \ref{allowedgamma.lemma}. 
With this choice of $N$ we have $\lambda_{\sharp} \asymp \sqrt {\log p / p } $. Recall that
according to Remark \ref{high-dim.remark} it must be true that $p \ge 1/ \lambda^{\sharp 2} $.
In the present context we in fact have $p \asymp \log p / \lambda^{\sharp 2} $.

\subsubsection{Creating $\gamma^0$ directly}

We first recall that 
$$\Theta_{1,1} - \Theta_{1,1}^{\sharp} \gg 0 \Leftrightarrow 
\| \Sigma_{-1,-1}^{1/2} (\gamma^0 - \gamma^{\sharp}) \|_2^2 \gg 0 $$
with $\Theta_1^{\sharp}$ given in (\ref{Thetasharp.equation}).
It will be the case in the constructions of this subsection. 

Fix some set $S\subset \{ 2 , \ldots , p \} $ with cardinality ${\rm s}$ and fix some $\lambda^{\sharp}$. We will assume 
$\lambda^{\sharp} \sqrt {\rm s} \rightarrow 0 $.

\begin{lemma} \label{direct.lemma}
Suppose there exists a vector $z \in \R^{p-1}$ with 
$$\| z \|_{\infty} \le 1 , 
$$
and 
\begin{equation}\label{z.equation}
1 - \lambda^{\sharp2} \| \Sigma_{-1,-1}^{-1/2} z \|_2^2 \gg 0  , \
 \lambda^{\sharp 2} \| \Sigma_{-1,-1}^{-1/2} z \|_2^2 \gg 0  . \
 \end{equation}
Let $\gamma^{\sharp}$ be a vector in $\R^{p-1}$ with 
$\gamma_{-S}^{\sharp}=0$ (i.e.\ $\gamma^{\sharp} = \gamma_S^{\sharp}$) and
$$ 1-  \lambda^{\sharp2} \| \Sigma_{-1,-1}^{-1/2} z \|_2^2 -\| \Sigma_{-1,-1}^{1/2}  \gamma_S^{\sharp} \|_2^2  \gg 0
 . $$
Define
$$\gamma^0 :=\gamma^{\sharp}+ \lambda^{\sharp}\Sigma_{-1,-1}^{-1} z . $$
Then, if $\lambda^{\sharp}\sqrt {\rm s} \rightarrow 0$, the pair $(\gamma^{\sharp}, \lambda^{\sharp})$ is an
eligible pair. Moreover, $\gamma^0$ is eventually allowed, $\lambda^{\sharp} \| \gamma^0 \|_1 \not\rightarrow 0 $
and in fact
$$ \| \Sigma_{-1,-1}^{1/2} (\gamma^0 - \gamma^{\sharp}) \|_2^2 \gg 0. $$
\end{lemma}

\begin{remark}\label{Lambdamax.remark}
Recall that $\Lambda_{\rm min}^2 (\Sigma_{-1,-1} ) \gg 0 $. To check condition (\ref{z.equation}) for
some $\| z \|_{\infty} \le 1$ one may want to impose in addition that
$\Lambda_{\rm max} (\Sigma_{-1,-1} ) = {\mathcal O} (1)$.  Then the condition is true
if $\| z \|_2  \asymp 1/\lambda^{\sharp}$ which is true for instance if
 $1/ \lambda^{\sharp2}$ of the coefficients of $z$ stay away from zero.  This is only possible
 if $p > 1/ \lambda^{\sharp 2}$ i.e., in high-dimensional situations (see also
 Remark \ref{high-dim.remark}). 
 \end{remark}
 
 \begin{remark} \label{random.remark} One may also consider taking
 $z= Z_{-1}^T \xi/(N \lambda^{\sharp})$ where $Z_{-1}$, $\xi $ and $\lambda^{\sharp}$ are
 is in Subsection \ref{regression.section}. Then
 $$ \lambda^{\sharp 2} \| \Sigma_{-1,-1}^{1/2} z \|_2^2 = \chi_{p-1}^2 / N $$
 as there and one arrives at eligible pairs ``with high probability". 
 \end{remark}
 
We now examine the following question: can one choose $\gamma_S^{\sharp}$ in Lemma
\ref{direct.lemma} equal to $\gamma_S^0$?
As we will see this will only be possible if a form of the irrepresentable condition
 holds.
The ``usual" irrepresentable condition (that implies
the absence of false positives of the Lasso, see \cite{ZY07}) involves the coefficients of the 
projection of the ``large" collection ${\bf x}_{-S} $ on the ``small" collection ${\bf x}_S$.
In our case, we reverse the roles of $S$ and $-S$.

\begin{definition} Fix a vector $z_{-S} \in \R^{p-{\rm s}-1}$ with $\| z_{-S} \|_{\infty} \le 1$. We say that the
reversed irrepresentable condition holds for $(S, z_{-S})$ if
$$\| \Sigma_{S, -S} \Sigma_{-S, -S}^{-1} z_{-S} \|_{\infty} \le 1 . $$
\end{definition}

\begin{lemma}\label{reversed-irrepresentable.lemma} $ $\\
a) Assume the reversed irrepresentable condition holds for $(S , z_{-S})$, 
and that in addition
$$1-   \lambda^{\sharp 2}\| \Sigma_{-1,-1}^{-1/2}  z_{-S} \|_2^2 \gg 0  , \ 
\lambda^{\sharp2} \| \Sigma_{-1,-1}^{-1/2}  z_{-S} \|_2^2 \gg 0  . $$
Let $\gamma_S^0 $ satisfy
$$1- \lambda^{\sharp2} \| \Sigma_{-1,-1}^{-1/2}  z_{-S} \|_2^2-
\| \Sigma_{-1,-1}^{1/2}  \gamma_S^0 \|_2^2 
\gg 0  . $$
Define $\gamma_{-S}^0 :=  \lambda^{\sharp}  \Sigma_{-S,-S}^{-1} z_{-S}$ and 
$\gamma^{\sharp}:= \gamma_S^0$
Then, if $\lambda^{\sharp} \sqrt {\rm s} \rightarrow 0$, the pair $(\gamma^{\sharp}, \lambda^{\sharp})$
is an eligible pair, $\gamma^0$ is eventually allowed and 
$\lambda^{\sharp} \| \gamma^0 \|_1 \not\rightarrow 0 $. In fact
$$ \| \Sigma_{-1,-1}^{1/2} (\gamma^0 - \gamma^{\sharp}) \|_2^2 \gg 0. $$
b) Conversely, if for some $\gamma^0$ and for $\gamma^{\sharp} := \gamma_S^0$ 
the pair $(\gamma^{\sharp}, \lambda^{\sharp})$
is an eligible pair, then the reversed irrepresentable condition holds for $(S , z_{-S})$ with appropriate
$z_{-S}$ satisfying $\| z_{-S} \|_{\infty} \le 1 $. 
\end{lemma}

\subsubsection{Creating $\gamma^0$ using a non-sparsity restriction} 

Fix some set $S\subset \{2 , \ldots , p \}$ with cardinality ${\rm s}$ and some $\lambda^{\sharp}>0$.
Let $w \in \R^{p-1}$ be a vector of strictly positive weights with $\| w \|_{\infty} \le 1$ and define
the matrix $W$ as the diagonal matrix with $w$ on its diagonal.
 Let
$$ c^0 \in \arg \min \biggl \{ \| \Sigma_{-1,-1}^{1/2} {c} \|_2^2 :\ \lambda^{\sharp} \| ( W {c} )_{-S} \|_1 = 1 \biggr \}  $$
and
$$\zeta_S := 0, \  \zeta_j  := {\rm sign} (c_j^0 ), \ j \notin S . $$

\begin{lemma}\label{Lagrangian.lemma}
The random variable
${\bf x}_{-1} c^0 $ is orthogonal to (i.e.\ independent of) ${\bf x}_S $.
Moreover
$$\| \Sigma_{-1,-1} c^0 \|_{\infty} = { \| w_{-S} \|_{\infty}\over  \lambda^{\sharp} \|  \Sigma_{-1,-1}^{-1/2} W \zeta \|_2^2 } $$
and
$$\| \Sigma_{-1,-1}^{1/2} c^0 \|_2^2 = {  1 \over \lambda^{\sharp2} \|  \Sigma_{-1,-1} ^{-1/2} W \zeta \|_2^2 } $$
\end{lemma}

\begin{lemma} \label{non-sparse.lemma} Suppose that
\begin{equation}\label{weights.equation}
1-  (\lambda^{\sharp2} \| \Sigma_{-1,-1}^{-1/2} W \zeta \|_2^2)^{-1} \gg 0   
 , \  \lambda^{\sharp2} \| \Sigma_{-1,-1}^{-1/2} W \zeta \|_2^2 = {\mathcal O} (1) . 
 \end{equation}
Let $\gamma^{\sharp}$ be a vector satisfying $\gamma_{-S}^{\sharp}=0$ and
$$0<1-   (\lambda^{\sharp2} \| \Sigma_{-1,-1}^{-1/2}W  \zeta \|_2^2)^{-1}-
\| \Sigma_{-1,-1}^{1/2}  \gamma_S^{\sharp} \|_2^2 \gg 0. $$
Define $\gamma_{-S}^0 = c_{-S}^0 $ and ${\bf x}_{-1} \gamma^0 :=
{\bf x}_S \gamma_S^{\sharp}  + ({\bf x}_{-S} {\rm A} {\bf x}_S ) \gamma_{-S}^0 $. 
Then ${\bf x}_S \gamma_S^{\sharp} $ is the projection of ${\bf x}_{-1} \gamma^0 $ on
${\bf x}_S$. Moreover, if $\lambda^{\sharp}\sqrt {\rm s} \rightarrow 0$, the pair
$(\gamma^{\sharp}, \lambda^{\sharp})$ is eligible, $\gamma^0$ is allowed
and $\lambda^{\sharp} \| \gamma^0 \|_1 \not\rightarrow 0 $.
In fact
$$ \| \Sigma_{-1,-1}^{1/2} (\gamma^0 - \gamma^{\sharp}) \|_2^2 \gg 0. $$
\end{lemma}

\begin{remark}
As in  Remark \ref{Lambdamax.remark}, 
to deal with the requirement  (\ref{weights.equation}) one may want to impose the condition 
$\Lambda_{\rm max}^2 ( \Sigma_{-1,-1} )={\mathcal O} (1)$.
\end{remark} 

\section{The case $\Sigma $ unknown}\label{unknownSigma.section}

As we will see
the concept of an eligible pair will also play a crucial role when
$\Sigma$ is unknown.

We use in this section the noisy Lasso
\begin{equation}\label{noisyLasso.equation}
\hat \gamma \in \arg \min_{c \in \R^{p-1} } \biggl \{
\| X_{-1} - X_{-1} c \|_2^2 / n + 2  \lambda^{{\rm Lasso}} \| c \|_1 \biggr \} . 
\end{equation}
We require the the tuning
parameter $\lambda^{{\rm Lasso}} $ to be 
of order $\sqrt {\log p / n}$.
Then in the debiased Lasso given in (\ref{de-biased.equation}) we apply
$$ \tilde \Theta_1 := \hat \Theta_1 $$ where
\begin{equation}\label{Thetahat.equation}
\hat  \Theta_1 := \begin{pmatrix} 1 \cr - \hat \gamma \end{pmatrix}/ (\| X_1 - X_{-1} \hat \gamma \|_2^2 / n +
 \lambda^{\rm Lasso} \| \hat  \gamma \|_1  ). 
 \end{equation}

Let $(\gamma^{\sharp}, \lambda^{\sharp})$ be an
eligible pair, with $\lambda^{\sharp} = {\mathcal O} ( \sqrt {\log p / n} ) $. 
In Lemma \ref{slowLasso.lemma}, we in fact need $\lambda^{\rm Lasso}$ to be at least as large as $\lambda^{\sharp}$. The
latter is allowed to be of small order $\sqrt {\log p / n}$, yet we will assume
$\lambda^{{\rm Lasso}} \| \gamma^{\sharp} \|_1 \rightarrow 0 $.

Let $\eta_n \rightarrow 0 $ be a sequence 
such that 
\begin{equation}\label{quadratic-forms.equation}
 \PP \biggl ( 
\inf_{c: \ \lambda^{{\rm Lasso}} \| c \|_1 \le 4\eta_n^2  , \ \| \Sigma_{-1,-1}^{1/2} c\|_2  = 1 } \| X_{-1} c \|_2^2/n  \ge
{1 \over 2}  \biggr ) \rightarrow 1 .
\end{equation}
Such a sequence exists, see for example Chapter 16 in \cite{vdG2016} and its references.

Define
 ${\varepsilon} := X_{1}- X_{-1} \gamma^0$ and
${\varepsilon}^{\sharp} := X_1- X_{-1} \gamma^{\sharp}= \varepsilon + X_{-1}
(\gamma^0 - \gamma^{\sharp} )$.  Choose $\lambda_{\varepsilon}^{\sharp} \asymp
\sqrt {\log p / n} $ in such a way that
$$\PP \biggl ( \| X^T \epsilon^{\sharp} \|_{\infty}/n\ge \lambda_{\varepsilon}^{\sharp}   \biggl ) \rightarrow 0$$
(see Lemma \ref{lambda*.lemma}). This implies $\lambda_{\epsilon}^{\sharp} \ge \lambda^{\sharp}$. 

The following lemma establishes the so-called ``slow rate". The proof is standard (see for example
\cite{BvdG2011}, Theorem 6.3), up to the replacement of $\varepsilon$ by
$\varepsilon^{\sharp}$. The result is the noisy counterpart of Lemma \ref{findpair.lemma}

\begin{lemma}\label{slowLasso.lemma} Let $\lambda^{{\rm Lasso}} \asymp \sqrt {\log p / n }$ satisfy
$\lambda^{{\rm Lasso}}  \ge 2 \lambda_{\varepsilon}^{\sharp} $ and suppose
$\lambda^{{\rm Lasso}}$ $\| \gamma^{\sharp} \|_1 \le \eta_n^2  $. Then
we have
$$ \lambda^{\rm Lasso } \| \hat \gamma \|_1 = o_{\PP}(1) , \  
 \| X_{-1}  ( \hat \gamma - \gamma^{\sharp}) \|_2^2/n= o_{\PP} (1) , \ 
\| \Sigma_{-1,-1}^{1/2}  ( \hat \gamma - \gamma^{\sharp}) \|_2^2 = o_{\PP} (1) . $$
\end{lemma}


\begin{theorem} \label{noisy.theorem}
Let $\hat b_1$ be the debaised Lasso given in (\ref{de-biased.equation})  with $\tilde \Theta_1$ equal to
$\hat \Theta_1$ given in
(\ref{Thetahat.equation}):
$$ \hat b_1 := \hat \beta_1 + \hat \Theta_1 X^T (Y- X \hat \beta )/n . $$
Assume that $\lambda^{\sharp} = {\mathcal O} (\sqrt {\log p / n})$, that
$\sqrt {\log p / n} \| \gamma^{\sharp} \|_1 = o(1)$ and 
$\lambda^{\rm Lasso} \asymp \sqrt {\log p / n} $ is sufficiently large 
(depending only on $\lambda^{\sharp}$).
Suppose that uniformly in $\beta^0 \in {\cal B}$
\begin{equation}\label{betarate.equation}
 \sqrt {\log p }  \| \hat \beta - \beta^0 \|_1= o_{\PP_{\beta^0}} (1).
 \end{equation}
 Then uniformly in $\beta^0 \in {\cal B}$
 $$ \hat b_1 - \beta_1^0 = \hat \Theta_1^T X^T \epsilon / n + o_{\PP_{\beta^0}} (1/\sqrt n) .$$
 Moreover,
$$\lim_{n \rightarrow\infty}
\sup_{\beta^0 \in {\cal B} } \PP \biggl ( \sqrt n {  (\hat b_1 - \beta_1^0 )\over\sqrt { \hat \Theta_1^T \hat \Sigma \hat \Theta_1 } } \le z \biggr )=  \Phi (z) \ \forall  \ z \in \R , $$
and for $\Theta_1^{\sharp} $ given in (\ref{Thetasharp.equation})
$$ \hat \Theta_1^T \hat \Sigma \hat \Theta_1= \Theta_{1,1}^{\sharp} + o_{\PP} (1) . $$
\end{theorem}

\begin{remark} Recall that Lemmas \ref{direct.lemma}, \ref{reversed-irrepresentable.lemma} and
\ref{non-sparse.lemma} present examples of eligible pairs
$(\gamma^{\sharp}, \lambda^{\sharp})$ where $\Theta_{1,1}^{\sharp}$ 
remains strictly smaller than $\Theta_{1,1}$. 
\end{remark}

\begin{remark} By Slutsky's Theorem we conclude that the asymptotic variance of $\hat b_1$
is (up to smaller order terms) equal to $\Theta_{1,1}^{\sharp}$. 
\end{remark}

\begin{remark} Note that Theorem \ref{noisy.theorem} only requires sparsity of $\gamma^{\sharp}$ in $\ell_1$-sense:
it requires $\| \gamma^{\sharp} \|_1= o( \sqrt {n / \log p})$. Indeed, we in fact base
the result on ``slow rates" for the Lasso $\hat \gamma$ given in Lemma \ref{slowLasso.lemma}
\end{remark}

\begin{remark}\label{asymptotic-linearity.remark} Assume the conditions of Theorem \ref{noisy.theorem}.
and that in fact
\begin{equation}\label{gammarate.equation}
\sqrt {\log p} \| \hat \gamma - \gamma^{\sharp} \|_1 = o_{\PP} (1) . 
\end{equation}
Then also $\sqrt {\log p} \| \hat \Theta_1  - \Theta_1^{\sharp} \|_1 = o_{\PP} (1)$, which implies
$$ \hat \Theta_1 X^T \epsilon /n = \Theta_1^{\sharp} X^T \epsilon / n + o_{\PP} (1/ \sqrt n) . $$
Thus, then the estimator $\hat b_1$ is asymptotically linear, uniformly in $\beta^0 \in {\cal B}$. 
The uniform asymptotic linearity of $\hat b_1$ implies in turn that the Cram\'er Rao lower bound of Subsection
\ref{CRLB.section} applies. 
\end{remark}

The results of the following two examples are summarized in Table \ref{conditionsunknown.table}.

\begin{example} \label{l0-unknown.example}
Assume the sparse model (\ref{sparse-model.equation}) with
$s = o ( \sqrt n / \log p ) $. As stated in Example \ref{l0.example}, for the 
Lasso estimator 
$\hat \beta $ given in (\ref{Lasso.equation}),
with appropriate choice of the tuning parameter $\lambda \asymp \sqrt {\log p / n } $, one has
uniformly in $\beta^0 \in {\cal B}$
$$\| \hat \beta - \beta^0 \|_1 = O_{\PP_{\beta^0}} ( s \sqrt {  \log p / n} )=
o_{\PP_{\beta^0}} (1/\sqrt {\log p} ) . $$
Suppose now as in Theorem \ref{noisy.theorem} that $\lambda^{\sharp} = {\mathcal O}
(\sqrt {\log p / n })$ and $\sqrt {\log p / n } \| \gamma^{\sharp} \|_1 = o(1)$. 
In fact, assume that $\| \gamma^{\sharp} \|_0^0$ is small enough so that
$\Theta_1^{\sharp}$ is a model direction.
Then one obtains
by the same arguments  if $\lambda^{\rm Lasso} \asymp
\sqrt {\log p /n }$ is suitably chosen
$$\| \hat \gamma - \gamma^{\sharp} \|_1 = O_{\PP} ( s \sqrt {  \log p / n}  )= o_{\PP_{\beta^0}}
(1/ \sqrt {\log p} ).$$
Thus then (\ref{gammarate.equation}) is met so that we have asymptotic linearity.
It means that the Cram\'er Rao lower bound applies and is achieved.
\end{example}

The model (\ref{l1-model.equation}) is too large for us to be able to apply when $\Sigma$ is unknown.
We now turn to the model (\ref{lr-model.equation}).

\begin{example}\label{lr-unknown.example}
For the model (\ref{lr-model.equation}) the rates for the Lasso $\hat \beta$ are given in Example
\ref{lr.example}. One sees that for $0 \le r < 1$ the requirement on $s$ becomes
$$ s = o(n^{1-r \over 2-r} / \log p ) . $$
By the same arguments, if for a fixed $0 \le {\rm r} \le 1 $
$$\| \gamma^{\sharp} \|_{\rm r}^{\rm r}= o(n^{1-{\rm r} \over 2} / \log^{2- {\rm r} \over 2} p )$$
one finds
$$ \| \hat \gamma - \gamma^{\sharp} \|_1 =  o_{\PP} (1/ \sqrt {\log p} ) $$
which yields asymptotic linearity so that the Cram\'er Rao lower bound applies. If such $\gamma^{\sharp}$ is in
addition a model direction after scaling, i.e. if
$\| \gamma^{\sharp} \|_r^r = {\mathcal O} ( n^{r \over 2} s^{2-r \over 2} )$, then
the Cramer Rao lower bound is achieved whenever $\beta^0$ stays away from the boundary. 
\end{example}

\section{Conclusion}\label{conclusion.section}
This paper illustrates that $\Theta_{1,1}$ can be larger than the asymptotic Cram\'er
Rao lower bound, that for certain $\Sigma$ the
asymptotic variance of a debiased Lasso is smaller than $\Theta_{1,1}$ and that in special such cases the asymptotic Cram\'er Rao lower bound
is achieved. 
In Examples \ref{l1-model.example} and \ref{lr-model.example} we showed that if $\beta^0$ stays
away from the boundary, then the asymptotic Cram\'er Rao lower bound is
$$\biggl (  \min_{\| c \|_r^r \le M n^{r \over 2} s^{2-r \over 2}}  \EE ( {\bf x}_1 - {\bf x}_{-1} c )^2 \biggr )^{-1}  $$
where $M$ is any fixed value. When $\Sigma$ is known, Theorem \ref{bsharp.theorem} shows
that up to log-terms this lower bound is achieved as soon as there exists an eligible pair $(\gamma^{\sharp}, \lambda^{\sharp})$.
with $\| \gamma^{\sharp} \|_r^r = {\mathcal O} (n^{r \over 2} s^{2-r \over 2} ) $.
When $\Sigma$ is unknown the situation is more involved and in particular for model (\ref{sparse-model.equation})
sparsity variant (i) is replaced by the stronger variant (ii). Model (\ref{l1-model.equation}) is too large for the case $\Sigma$ unknown and model (\ref{lr-model.equation}) requires more sparsity then
model (\ref{sparse-model.equation}): the larger $r$ the smaller $s$ is required to be. 
Model (\ref{sparse-model.equation}) however appears for both known and unknown $\Sigma$ too
stringent as results depend on the exact value of $s$, not only on its order.
Model (\ref{l1-model.equation}) (with $\Sigma $ known) or more generally model (\ref{lr-model.equation})
(with $0<r<1$ if $\Sigma$ is unknown)
do not suffer from such a dependence as long as $\beta^0$
stays away from the boundary.  

\section{Proofs} \label{proofs.section}

\subsection{Proof for Section \ref{introduction.section}} 

The proof of Proposition  \ref{CRLB.proposition}
relies on the results in \cite{Jankova16}, which allow the arguments to
follow those of the
low-dimensional case. These arguments are then rather standard.

{\bf Proof of Proposition \ref{CRLB.proposition}.}  Let $h \in {\cal H}_{\beta^0}$, $|h_1 | \ge \rho$ and
$h^T \Sigma h \le R^2 $.
The log-likelihood ratio ${\cal L} (h)$ for  $\beta^0+ h / \sqrt n$
with respect to $\beta^0$ is
\begin{eqnarray*}
{\cal L} (h) &=& h^T X^T \epsilon / \sqrt n + h^T \hat \Sigma h / 2\\
& = & h^T X^T \epsilon / \sqrt n + h^T \Sigma h / 2+ o_{\PP} (1) 
\end{eqnarray*}
since $h^T \Sigma h = {\mathcal O} (1)$. 
Let for $(x , y) \in \R^{p+1} $
$$ l_{\beta_0} (x,y):= \begin{pmatrix}{\bf i}_{\beta^0} (x,y) \cr x^T (y- x\beta^0) \end{pmatrix} $$
and define
$$\Omega_{\beta^0} := \EE_{\beta^0} l_{\beta^0} ({\bf x}, {\bf y} ) l_{\beta^0} ({\bf x}, {\bf y} )^T.$$ 
By the Lindeberg condition, we can apply Lindeberg's central limit theorem to conclude that
for any sequence $a:= (1, c^T  )^T \in \R^{p+1}$ with $c^T \Sigma c  = {\mathcal O} (1)$ it holds
that
$${a^T\sum_{i=1}^n  l_{\beta^0} (X_i, Y_i) \over \sqrt { n a^T  \Omega_{\beta^0} 
a}} ){\buildrel{{\cal D}_{ \beta^0    }}\over \longrightarrow} {\cal N} (0,1)  .
$$
Therefore, by Wold's device
$$\biggl (  \begin{pmatrix} 1 & 0 \cr 0 & h^T\cr \end{pmatrix}\Omega_{\beta^0}  \begin{pmatrix} 1 & 0 \cr 0 & h
\end{pmatrix} \biggr )^{-1/2} 
 \begin{pmatrix} 1 & 0 \cr 0 & h^T\cr \end{pmatrix} \sum_{i=1}^n l_{\beta^0 } (X_i, Y_i) / 
 \sqrt n\ \ {\buildrel{{\cal D}_{ \beta^0    }}\over \longrightarrow} \ {\cal N} \biggl (0,
 \begin{pmatrix}1 & 0 \cr 0 & 1 \cr  \end{pmatrix} \biggr ) .  $$
 We now apply  a slight modification of Lemmas 16 and 23 in \cite{Jankova16}, where we
 drop the assumption of bounded eigenvalues of $\Sigma$ 
 (which is possible because we have $h^T \Sigma h \le R^2 = {\mathcal O} (1)$). The asymptotic linearity of $T$ 
 and the 2-dimensional central limit theorem just obtained imply that
 at the alternative $\beta^0 + h / \sqrt n $ it holds that 
 $$ { T - (\beta_1^0 + h_1 / \sqrt n) +h_1  - h^T\EE_{\beta^0} {\bf i}_{\beta^0}  ({\bf x}, {\bf y} )
  {\bf x}^T ({\bf y } - {\bf x} \beta^0 ) \over \sqrt {n V_{\beta^0}^2 } } \ \ 
 {\buildrel{{\cal D}_{ \beta^0  + h / \sqrt n   }}\over \longrightarrow} {\cal N} (0 ,1) . $$
 As $T$ is assumed to be regular at $\beta^0$ we conclude that
 $$h^T \EE_{\beta^0} {\bf i}_{\beta^0}  ({\bf x}, {\bf y} )
  {\bf x}^T ({\bf y } - {\bf x} \beta^0 ) = h_1 + o(1) . $$
  But by the Cauchy-Schwarz inequality
 $$ \biggl (  h^T \EE_{\beta^0} {\bf i}_{\beta^0}  ({\bf x}, {\bf y} )
  {\bf x}^T ({\bf y } - {\bf x} \beta^0 ) \biggr )^2 \le V_{\beta^0}^2 h^T \Sigma h .$$
  Moreover, 
    $$\biggl (h_1 +o(1)\biggr )^2= h_1^2 + o(1) . $$
  so that we obtain
  \begin{equation}\label{oneh.equation}
  V_{\beta^0}^2  \ge { h_1^2 +o(1) \over h^T \Sigma h } = { h_1^2 \over h^T \Sigma h} + o(1)   
  \end{equation}
  where in the last step we used $h^T \Sigma h \ge \| h\|_2^2 / \Lambda_{\rm min}^2 \ge \rho^2$
  so that  $1/h^T \Sigma h= {\mathcal O} (1)$. 
   Since the result is true for all $h \in {\cal H}_{\beta^0}$, $|h_1 | \ge \rho$ and
$h^T \Sigma h \le R^2 $, we may maximize the right hand side of
(\ref{oneh.equation}) over all such $h$.
 \hfill $\sqcup \mkern -12mu \sqcap$

\subsection{Proofs for Section \ref{knownSigma.section}} 

 {\bf Proof of Lemma \ref{need.lemma}.} 
  Because ${\bf x}_{-1} \gamma^0$ is the projection of ${\bf x}_1 $ on
  ${\bf x}_{-1} $ we know that
  $$\EE ( {\bf x}_1 - {\bf x}_{-1} \gamma^{\sharp} )^2 \ge
  \EE ( {\bf x}_1 - {\bf x}_{-1} \gamma^0 )^2 = 1/ \Theta_{1,1} . $$
  Moreover
  $$\Theta_{1,1} = {\rm e}_1^T \Theta {\rm e}_1 \le \Lambda_{\rm max}^2 (\Theta ) =
  1/ \Lambda_{\rm min}^2 . $$
  Thus
  \begin{equation}\label{projectx1.equation}
  \EE ( {\bf x}_1 - {\bf x}_{-1} \gamma^{\sharp} )^2  \ge \Lambda_{\rm min}^2 . 
  \end{equation}
  We now rewrite
  \begin{eqnarray*}
  \EE ( {\bf x}_1 - {\bf x}_{-1} \gamma^{\sharp} )^2 &=& 1+
  \gamma^{\sharp T} \Sigma_{-1,-1} \gamma^{\sharp} - 2 \EE {\bf x}_1^T {\bf x}_{-1} \gamma^{\sharp} \\
  &=& 1+  \gamma^{\sharp T} \Sigma_{-1,-1} \gamma^{\sharp} - 2 \gamma^{0T}
  \Sigma_{-1,-1} \gamma^{\sharp} 
  \end{eqnarray*}
  where in the second equality we
  used that ${\bf x}_{1} - {\bf x}_{-1} \gamma^0$ is the anti-projection
  of ${\bf x}_1 $ on ${\bf x}_{-1}$  and hence orthogonal to ${\bf x}_{-1} \gamma^{\sharp}$. 
  For the cross-product we have by the two conditions on the pair $(\gamma^{\sharp}, \lambda^{\sharp})$
  \begin{eqnarray*}
   \gamma^{0T}
  \Sigma_{-1,-1} \gamma^{\sharp} &= &\gamma^{\sharp T} \Sigma_{-1, -1} \gamma^{\sharp} +
  \underbrace{( \gamma^0 - \gamma^{\sharp} )^T \Sigma_{-1, -1} \gamma^{\sharp} }_{ | \cdot | \le
 \|  \Sigma_{-1, -1}( \gamma^0 - \gamma^{\sharp} ) \|_{\infty} \|  \gamma^{\sharp}\|_1 \le
 \lambda^{\sharp} \| \gamma^{\sharp} \|_1= o(1) } \\
 &=& \gamma^{\sharp T} \Sigma_{-1, -1} \gamma^{\sharp} + o(1) .
  \end{eqnarray*}
  Thus
  $$\EE ( {\bf x}_1 - {\bf x}_{-1} \gamma^{\sharp} )^2 = 1- \gamma^{0 T} \Sigma_{-1, -1} \gamma^{\sharp}  +o(1) . $$
  Combining this with inequality (\ref{projectx1.equation}) proves the first result of the lemma. 
  The second result:
  $$ \| \Sigma( \Theta_1^{\sharp}- \Theta_1 )\|_{\infty} \le \lambda_0^{\sharp} $$
  follows trivially from this. 
  For the third result, we compute and re-use the already obtained results:
    \begin{eqnarray*}
   \Theta_1^{\sharp} \Sigma\Theta_1^{\sharp}& =&
  {1 - \gamma^{0 T} \Sigma_{-1,-1} \gamma^{\sharp} - \gamma^{\sharp T} \Sigma_{-1,-1} (\gamma^0 - \gamma^{\sharp} )\over
  (1 - \gamma^{0 T} \Sigma_{-1,-1} \gamma^{\sharp})^2}\\ & =&
   \underbrace{{1  \over
  1 - \gamma^{0 T} \Sigma_{-1,-1} \gamma^{\sharp}}}_{= \Theta_{1,1}^{\sharp} }  +o(1) \\
  &=& 1/ \EE ({\bf x}_1- {\bf x}_{-1} \gamma^{\sharp} )^2 + o(1) \\
  & \le & \underbrace{1/ \EE ({\bf x}_1- {\bf x}_{-1} \gamma^0 )^2 }_{= \Theta_{1,1} } +o(1) .
   \end{eqnarray*} 
   To show the final statement of the lemma, assume on the contrary
   that $\gamma^0$ is sparse:
   $$ \lambda^{\sharp} \| \gamma^0 \|_1 \rightarrow 0 . $$
   Then
   \begin{eqnarray*}
   \Theta_{1,1}^{\sharp} &= &{1 \over 1- \gamma^{0T} \Sigma_{-1,-1} \gamma^{\sharp} } \\
   & = & {1 \over 1- \gamma^{0T} \Sigma_{-1,-1} \gamma^0 + \underbrace{\gamma^{0T} \Sigma_{-1,-1}(\gamma^0 -
   \gamma^{\sharp} }_{|\cdot | \le \lambda^{\sharp } \| \gamma^0 \|_1 = o(1)} )}\\
   &=& \Theta_{1,1} + o(1) .
   \end{eqnarray*}
   \hfill $\sqcup \mkern -12mu \sqcap$
   
   {\bf Proof of Theorem \ref{bsharp.theorem}.} We use the decomposition
of the beginning of this section applied to $\hat b_{I,1}^{\sharp}$
$$\hat b_{I,1}^{\sharp} - \beta_1^0 =2 \Theta_1^{\sharp T} X_I^T \epsilon_I / n +  \underbrace{({\rm e}_1^T -
\Theta_1^{\sharp T} \hat \Sigma _I) ( \hat \beta_{II} - \beta_1^0 )}_{{\rm remainder}_I}  $$
where ${\rm remainder}_I$  is
$$ ({\rm e}_1^T-\Theta_1^{\sharp T}\hat \Sigma_I ) ( \hat \beta_{II} - \beta^0 )=
 \underbrace{ \Theta_1^{\sharp T} ( \Sigma - \hat \Sigma_I )(\hat \beta_{II} - \beta^0 )}_{:=(i)} + 
\underbrace{( {\rm e}_1^T  -  \Theta_1^{\sharp T}  \Sigma) (\hat \beta_{II} - \beta^0 )}_{:=(ii)} . 
 $$
 Here, $\epsilon_I := Y_I - X_I \beta^0$ and $\hat \Sigma_{I}:= 2X_I^T X_{I} /n $.
 But, given $(X_{II} , Y_{II})$,
 $$ \hat\Theta_1^{\sharp T} \hat \Sigma_I (\hat \beta_{II} - \beta^0) =
 2 \Theta_1^{\sharp T} X_I^T X_I(\hat \beta_{II} - \beta^0) /n$$
 is the average of $n/2$ i.i.d.\ random variables which are the product of a random variable with the
 ${\cal N} (0, \Theta_1^{\sharp T} \Sigma  \Theta_1^{\sharp})$-distribution 
 and a ${\cal N} (0, \| \Sigma_{1/2} ( \hat \beta_{II} - \beta^0) \|_2^2)$-distributed random variable. 
 Since the variances satisfy $\Theta_1^{\sharp T} \Sigma \Theta_1^{\sharp}= \Theta_{1,1}^{\sharp} +
 o(1) = {\mathcal O} (1) $ and $\| \Sigma_{1/2} ( \hat \beta_{II} - \beta^0) \|_2^2=o_{\PP_{\beta^0}} (1)$\
 uniformly in $\beta^0 \in {\cal B}$ we see that
 $$ (i) = o_{\PP_{\beta^0}} (1/ \sqrt n ) $$
 uniformly in $\beta^0 \in {\cal B}$.  For the term $(ii)$ we use that
 $$ \| \Sigma \Theta_1^{\sharp} - {\rm e}_1\|_{\infty} = {\mathcal O} (\lambda^{\sharp} ) $$
 and the assumption
 $\lambda^{\sharp} \| \hat \beta_{II} - \beta^0 \|_1 =o_{\PP_{\beta^0 }} (1/ \sqrt n ) $ uniformly
 in $\beta^0 \in {\cal B}$. 
 This gives that uniformly
 in $\beta^0 \in {\cal B}$
$$ \hat b_{I,1}^{\sharp} - \beta_1^0 = 2 \Theta_1^{\sharp T}X_I^T \epsilon_I /n +
 o_{\PP_{\beta^0}} (1/ \sqrt n). $$
 In the same way one derives that
 uniformly
 in $\beta^0 \in {\cal B}$
$$ \hat b_{II,1}^{\sharp} - \beta_1^0 = 2 \Theta_1^{\sharp T} X_{II}^T \epsilon_{II}/n +
 o_{\PP_{\beta^0}} (1/ \sqrt n)$$
 with $\epsilon_{II} := Y_{II} - X_{II} \beta^0 $. Since
 $\hat b_1^{\sharp} = (b_{I,1}^{\sharp} + b_{II,1}^{\sharp} )/2$ is the average of the
 two, this proves the asymptotic linearity. Further
 ${\rm var} ( \Theta_1^{\sharp} X^T \epsilon / \sqrt n) = \Theta_1^{\sharp T} \Sigma \Theta_1^{\sharp}=
 \Theta_{1,1}^{\sharp}+o(1)$ by Lemma \ref{need.lemma}. The central limit theorem completes the proof. 
 \hfill $\sqcup \mkern -12mu \sqcap$
 
 \subsection{Proofs for Section \ref{eligible.section}}
 
 {\bf Proof of Lemma \ref{unique.lemma}.} We have
\begin{eqnarray*}
(\gamma^{\sharp} - \gamma^{\flat} )^T\Sigma_{-1,-1} (\gamma^{\sharp} - \gamma^{\flat} )
&\le & \|\gamma^{\sharp} - \gamma^{\flat} \|_1 \| \Sigma_{-1,-1} (\gamma^{\sharp} - \gamma^{\flat} )\|_{\infty} \\
& \le &\lambda^{\sharp} \|\gamma^{\sharp} - \gamma^{\flat} \|_1 \\
& \le & \lambda^{\sharp} \|\gamma^{\sharp}\|_1 +  \lambda^{\sharp} \|  \gamma^{\flat} \|_1 \rightarrow 0 
\end{eqnarray*}
\hfill $\sqcup \mkern -12mu \sqcap$

{\bf Proof of Lemma \ref{findpair.lemma}.}  By the KKT conditions
$$(\gamma_{\rm Lasso} - \gamma^{\sharp})^T \Sigma_{-1, -1} ( \gamma_{\rm Lasso} - \gamma^0) \le
\lambda_{\rm Lasso} \| \gamma^{\sharp} \|_1-  \lambda_{\rm Lasso} \| \gamma_{\rm Lasso} \|_1 .$$
Therefore,
\begin{eqnarray*}
& &  (\gamma_{\rm Lasso} - \gamma^{\sharp} )^T \Sigma_{-1,-1} (\gamma_{\rm Lasso} - \gamma^{\sharp} ) \\& =& 
 (  \gamma_{\rm Lasso} - \gamma^{\sharp})^T \Sigma_{-1, -1} ( \ \gamma_{\rm Lasso} - \gamma^0)+
 (  \gamma_{\rm Lasso} - \gamma^{\sharp})^T \Sigma_{-1, -1} ( \gamma^0 - \gamma^{\sharp})\\
 &\le & \lambda_{\rm Lasso} \| \gamma^{\sharp} \|_1-  \lambda_{\rm Lasso} \|  \gamma_{\rm Lasso} \|_1 +
  \|  \gamma_{\rm Lasso} - \gamma^{\sharp}\|_1 \|  \Sigma_{-1, -1} ( \gamma^0 - \gamma^{\sharp})\|_{\infty} \\
  & \le & \lambda_{\rm Lasso} \| \gamma^{\sharp} \|_1-  \lambda_{\rm Lasso} \| \gamma_{\rm Lasso} \|_1 +
  \lambda^{\sharp}  \|  \gamma_{\rm Lasso} - \gamma^{\sharp}\|_1  \\
  &\le& ( \lambda_{\rm Lasso} + \lambda^{\sharp} )\| \gamma^{\sharp} \|_1 -
  ( \lambda_{\rm Lasso} - \lambda^{\sharp} ) \|  \gamma_{\rm Lasso} \|_1 .
  \end{eqnarray*}
  Thus
  $$  (\gamma_{\rm Lasso} - \gamma^{\sharp} )^T\Sigma_{-1,-1}(\gamma_{\rm Lasso} - \gamma^{\sharp} ) +
   ( \lambda_{\rm Lasso} - \lambda^{\sharp} ) \| \gamma_{\rm Lasso} \|_1\le  ( \lambda_{\rm Lasso} + \lambda^{\sharp} )\| \gamma^{\sharp} \|_1 \rightarrow 0 $$
   where we used that $\lambda_{\rm Lasso} > \lambda^{\sharp}$ and $\lambda_{\rm Lasso} \| \gamma^{\sharp} \|_1 \rightarrow 0$.
      We also know that by the KKT conditions
   $$\|  \Sigma_{-1, -1} ( \gamma_{\rm Lasso}- \gamma^0) \|_{\infty} \le  \lambda_{\rm Lasso} . $$
   If
   $\lambda_{\rm Lasso} \ge 2 \lambda^{\sharp} $, we obtain from the above
   $$ \lambda_{\rm Lasso} \|  \gamma_{\rm Lasso} \|_1 \le 
   3 \lambda_{\rm Lasso}  \| \gamma^{\sharp} \|_1  /2  \rightarrow 0 . 
   $$
   So $( \gamma_{\rm Lasso}, \lambda_{\rm Lasso})$ is an eligible pair.
   \hfill $\sqcup \mkern -12mu \sqcap$
   
   \vskip .1in
   {\bf Proof of Lemma \ref{antiprojection.lemma}.} 
Note that
\begin{eqnarray*}
 \EE {\bf x}_{-S}^T {\bf x}_{-1} \gamma^{-S} &=&
\EE {\bf x}_{-S}^T \biggl [ ( {\bf x}_{-S} \gamma_{-S}^0 ) {\rm A} {\bf x}_S \biggr ]\\
&=& \EE {\bf x}_{-S}^T \biggl [ {\bf x}_{-S} - {\bf x}_S \Sigma_{S,S}^{-1} \Sigma_{S,-S}   \biggr ] \gamma_{-S}^0 \\
& = & \biggl [ \Sigma_{-S, -S} - \Sigma_{-S, S} \Sigma_{S,S}^{-1} \Sigma_{S,-S} \biggr ] \gamma_{-S}^0 . 
\end{eqnarray*} 
We therefore have
\begin{eqnarray*}
 \| v_{-S}^S \|_{\infty} &=& \biggl \| \biggl [ \Sigma_{-S, -S} - \Sigma_{-S, S} \Sigma_{S,S}^{-1} \Sigma_{S,-S} \biggr ] 
\gamma_{-S}^0 \biggr \|_{\infty} \\
& \le & \vvert \Sigma_{-S, -S} - \Sigma_{-S, S} \Sigma_{S,S}^{-1} \Sigma_{S,-S}\vvert_1 
\|\gamma_{-S}^0 \|_{\infty} . 
\end{eqnarray*}
\hfill $\sqcup \mkern -12mu \sqcap$

{\bf Proof of Lemma \ref{non-sparse1.lemma}.} 
We have 
$$\| \gamma_{-S}^0 \|_{\infty} = \lambda^S/C , $$
so that by Lemma \ref{antiprojection.lemma}
$$\| \Sigma_{-1,-1} \gamma^S \|_{\infty} \le \lambda^S . $$
Moreover
$$ \lambda^S \| \gamma^S \|_1 \le \lambda^S \sqrt {\rm s} \rightarrow 0 . $$
Thus $( \gamma^S , \lambda^S)$ is an eligible pair.
Finally
$$\lambda^S \| \gamma^0 \|_1 \ge
\lambda^S \| \gamma_{-S}^0 \|_1 = {\kappa (\gamma_{-S}^0) \over C   } \not\rightarrow 0 . $$
\hfill $\sqcup \mkern -12mu \sqcap$

{\bf Proof of Lemma \ref{projectS.lemma}.}
  Write
  $$ {\bf x}_{-1} \gamma^0 = {\bf x}_{-1} \gamma^{\sharp} + \xi^0 , \ \| \EE {\bf x}_{-1}^T \xi^0 \|_{\infty} \le \lambda^{\sharp}. $$
  It holds that
  $$ \gamma_S^S = \Sigma_{S,S}^{-1}  \EE  {\bf x}_S^T {\bf x}_{-1} \gamma^0 .$$
  So
  \begin{eqnarray*}
  \| \Sigma_{-1,-1}^{1/2} ( \gamma_S^S - \gamma^{\sharp} ) \|_2^2 &=&
  (\EE {\bf x}_S^T \xi^0  )^T\Sigma_{S,S}^{-1} ( \EE {\bf x}_S^T \xi^0)  \\
  & \le &  \| \EE {\bf x}_S ^T\xi^0 \|_2 / \Lambda_{\rm min}^2 \\
  & \le & {\rm s} \| \EE {\bf x}_S^T\xi^0  \|_{\infty}^2 / \Lambda_{\rm min}^2 \rightarrow 0 .
  \end{eqnarray*}
  \hfill $\sqcup \mkern -12mu \sqcap$ 

{\bf Proof of Lemma \ref{allowedgamma.lemma}.}

For $a \in \R$ and $c \in \R^{p-1}$ satisfying $a^2 + \| \Sigma_{-1,-1}^{1/2} c \|_2^2=1$,
\begin{eqnarray*}
 {\begin{pmatrix}  a\cr c \end{pmatrix}}^T \Sigma (\gamma^0) {\begin{pmatrix}  a\cr c \end{pmatrix}} 
& = & a^2+ 2 a \gamma^{0T} \Sigma_{-1,-1} c + c^T \Sigma_{-1,-1} c \\
& = &1+ 2 a \gamma^{0T} \Sigma_{-1,-1} c \\  &\ge & 1- 2|a|\| \Sigma_{-1, -1}^{1/2}\gamma^0 \|_2   \| \Sigma_{-1,-1}^{1/2} c \|_2 \\
&=& 1- 2 |a| \sqrt {1- a^2} \| \Sigma_{-1,-1}^{1/2} \gamma^0 \|_2 \\
&\ge & 1- \| \Sigma_{-1,-1}^{1/2} \gamma^0 \|_2. 
\end{eqnarray*} 

But then
\begin{eqnarray*}
 & \ & {1 \over a^2 + \| c \|_2^2 }  {\begin{pmatrix}  a\cr c \end{pmatrix}}^T \Sigma (\gamma^0) {\begin{pmatrix}  a\cr c \end{pmatrix}} \\
 &=&{a^2 +  \| \Sigma_{-1,-1}^{1/2} c \|_2^2 \over a^2 + \| c \|_2^2 }  {\begin{pmatrix}  a\cr c \end{pmatrix}}^T \Sigma (\gamma^0) {\begin{pmatrix}  a\cr c \end{pmatrix}} \\
 &\ge & { a^2 + \Lambda_{\rm min}^2 (\Sigma_{-1,-1} ) \| c \|_2 \over 
 a^2 + \| c \|_2 } (1- \| \Sigma_{-1,-1}^{1/2} \gamma^0 \|_2) \\
 &\ge & (1- \| \Sigma_{-1,-1}^{1/2} \gamma^0 \|_2) \Lambda_{\rm min}^2 (\Sigma_{-1,-1} ) .
\end{eqnarray*}

Hence $\Lambda_{\rm min}^2 (\Sigma (\gamma^0))$ is positive definite and 
$$\Lambda_{\rm min}^2 (\Sigma (\gamma^0)) \ge (1- \| \Sigma_{-1,-1}^{1/2} \gamma^0 \|_2 ) \Lambda_{\rm min}^2  
(\Sigma_{-1,-1} ) . $$
It further holds for all $j\in \{ 2 , \ldots , p \} $ that
$$ \EE| {\bf x}_j^T {\bf x}_{-1} \gamma^0 |\le \sqrt {\EE ({\bf x}_{-1} \gamma^0 )^2} < 1 . $$

\hfill $\sqcup \mkern -12mu \sqcap$

{\bf Proof of Lemma \ref{direct.lemma}.}
By definition
$$ \Sigma_{-1,-1} ( \gamma^0 - \gamma^{\sharp} ) = \lambda^{\sharp} z , $$
so that
$$ \| \Sigma_{-1,-1} ( \gamma^0 - \gamma^{\sharp} ) \|_{\infty} \le \lambda^{\sharp} .$$
Moreover, $\lambda^{\sharp} \| \gamma^{\sharp} \|_1 \le \lambda \sqrt{\rm s}  \| \gamma^{\sharp} \|_2 \rightarrow 0 $,
since 
$$\| \gamma^{\sharp} \|_2 \le \| \Sigma_{-1,-1}^{1/2} \gamma^{\sharp} \|_2 / \Lambda_{\rm min}
(\Sigma_{-1,-1} ) = {\mathcal O} (1) . $$
So $( \gamma^{\sharp}, \lambda^{\sharp})$ is an eligible pair.
We further have
$$ (\gamma^0 - \gamma^{\sharp})^T \Sigma_{-1,-1} (\gamma^0 - \gamma^{\sharp}) = \lambda^{\sharp2} z^T \Sigma_{-1, -1}^{-1} z . $$
Thus
\begin{eqnarray*}
\gamma^{0T} \Sigma_{-1,-1} \gamma^0& =& \gamma^{\sharp T} \Sigma_{-1-1} \gamma^{\sharp}+
2 ( \gamma^0 - \gamma^{\sharp})^T \Sigma_{-1,-1} \gamma^{\sharp} +
\lambda^{\sharp2} z^T \Sigma_{-1,-1}^{-1} z \\
&<& 1+ o(1)
\end{eqnarray*}
and 
$$1- \| \Sigma_{-1, -1}^{1/2} \gamma^0 \|_2^2 \gg 0  $$
where the positivity is true for large enough $n$. 
Therefore, by Lemma \ref{allowedgamma.lemma},  $\gamma^0$ is eventually allowed.
Finally, 
$$\lambda^{\sharp} \| \gamma^0 - \gamma^{\sharp}  \|_1 \ge \lambda^{\sharp} z^T ( \gamma^0 - \gamma^{\sharp})
= \lambda^{\sharp2 }z^T \Sigma_{-1,-1}^{-1} z \gg 0 . $$
So, since $\lambda^{\sharp} \| \gamma^{\sharp} \|_1 \rightarrow 0$, it must be true that
$\lambda^{\sharp} \| \gamma^0 \|_1 \gg 0 $.
We in fact have
$$\gamma^{0T} \Sigma_{-1,-1} \gamma^0 -
\gamma^{\sharp T} \Sigma_{-1,-1} \gamma^{\sharp} = \lambda^{\sharp2} z^T \Sigma_{-1, -1}^{-1} z + o(1)
 $$
 so that
$$  \gamma^{0T} \Sigma_{-1,-1} \gamma^0 -
\gamma^{\sharp T} \Sigma_{-1,-1} \gamma^{\sharp}  \gg 0 . $$
 
\hfill $\sqcup \mkern -12mu \sqcap$

{\bf Proof of Lemma \ref{reversed-irrepresentable.lemma}.} 

(${{\rm (a) } \atop \Rightarrow} $) 
For
$z_S :=  \Sigma_{S,-S} \Sigma_{-S,-S}^{-1} z_{-S}$
the equality 
$$ \Sigma_{-1,-1} (\gamma^0 - \gamma^{\sharp}) = \lambda^{\sharp} z  $$
holds. By assumption $\| z_{-S} \|_{\infty} \le 1 $ and by the 
reversed irrepresentable condition also
$ \| z_S \|_{\infty} \le 1 $. Thus
$$ \| \Sigma_{-1,-1}(\gamma^0 - \gamma^{\sharp})   \|_{\infty} \le \lambda^{\sharp} . $$

Moreover, 
$$\lambda^{\sharp} \| \gamma^{\sharp} \|_1 = 
\lambda^{\sharp} \| \gamma_S^0 \|_1 \le \lambda^{\sharp} \sqrt {\rm s} \| \gamma_S^0 \|_2 \rightarrow 0 . $$
 So $( \gamma^{\sharp} , \lambda^{\sharp}) $ is eligible.

To see make sure that $\gamma^0$ is allowed we bound $\gamma^{0T} \Sigma_{-1, -1}
\gamma^0$:
$$ \gamma^{0T} \Sigma_{-1,-1} \gamma^0 \le
\gamma_S^{0T} \Sigma_{S,S} \gamma_S^0 + 2 \lambda^{\sharp} \| \gamma_S^0 \|_1 +
\lambda^{\sharp 2} z_{-S}^T  \Sigma_{-1,-1}^{-1} z_{-S}. $$
Therefore, since $\lambda^{\sharp} \| \gamma_S^0 \|_1\rightarrow 0$,
and in view of 
Lemma \ref{allowedgamma.lemma}, the vector $\gamma^0$ is for large enough $n$ allowed. 
Finally, we have
\begin{eqnarray*}
\lambda^{\sharp} \| \gamma_0 \|_1& \ge&  \lambda^{\sharp} \| \gamma_{-S}^0 \|_1\\
& \ge&
\lambda^{\sharp}  z_{-S}^T \gamma_{-S}^0\\ & = & \lambda^{\sharp 2} z_{-S}^T \Sigma_{-S, -S}^{-1} z_{-S} 
\gg 0. 
\end{eqnarray*}
In fact
$$ \gamma^{0T} \Sigma_{-1,-1} \gamma^0 -
\gamma^{\sharp T} \Sigma_{-1,-1} \gamma^{\sharp} = 
 \lambda^{\sharp 2} z_{-S}^T \Sigma_{-S, -S}^{-1} z_{-S}+ o(1) \gg 0 
 . $$

(${ {\rm (b)} \atop \Leftarrow}$)
If $(\gamma^{\sharp}, \lambda^{\sharp})$ is an eligible pair we have
$$\| \Sigma_{-1,-1} ( \gamma^{\sharp} - \gamma^0 ) \|_{\infty} \le \lambda^{\sharp}. $$
Define now $ c =  \gamma^0 - \gamma^{\sharp} $ and
$z=  \Sigma_{-1,-1}^{-1} c / \lambda^{\sharp} $. 
Then 
$$ \Sigma_{-1,-1} c = \lambda^{\sharp} z $$
 and $\|z \|_{\infty} \le 1$, $c_S=0$.
It follows that
$$\Sigma_{S,-S} \Sigma_{-S,-S}^{-1} z_{-S} = z_S $$
so that
$$\| \Sigma_{S,-S} \Sigma_{-S,-S}^{-1} z_{-S} \|_{\infty} \le 1. $$
\hfill $\sqcup \mkern -12mu \sqcap$

{\bf Proof of Lemma \ref{Lagrangian.lemma}.} 
One readily verifies that
all $c_j^0$ with $j \notin S $ are non-zero.  One thus has
the Lagrangrian
$$\Sigma_{-1, -1}  c^0 = \tilde \lambda W \zeta   $$
where $\tilde \lambda $ is the Lagrangian parameter.
Since $\zeta_S =0$ this says that
$$(\Sigma_{-1,-1} c^0)_S = 0$$
so we know that ${\bf x}_{-1} c^0 $ is orthogonal to ${\bf x}_S$. 

The restriction gives
$$ \| \Sigma_{-1,-1}^{1/2} c^0 \|_2^2 = \tilde \lambda \| (W c^0)_{-S} \|_{\infty} = \tilde \lambda / \lambda^{\sharp} , $$
so
$$ \tilde \lambda= { \lambda^{\sharp} \| \Sigma_{-1,-1}^{1/2} c^0 \|_2^2 } , $$
and inserting this back yields
$$\Sigma_{-1,-1} c^0 =  \lambda^{\sharp} { \| \Sigma_{-1,-1}^{1/2} c^0 \|_2^2 }W \zeta. $$
It follows that
$$c^0 =\lambda^{\sharp} { \| \Sigma_{-1,-1}^{1/2} c^{0} \|_2^2 } \Sigma_{-1,-1}^{-1} W\zeta . $$
But then
$$ \| \Sigma_{-1,-1}^{1/2} c^0 \|_2 =\lambda^{\sharp} { \| \Sigma_{-1,-1}^{1/2} c^0 \|_2^2 } 
\| \Sigma_{-1,-1}^{-1/2} W \zeta\|_2  $$
or
$$\| \Sigma_{-1,-1}^{1/2} c^0\|_2 = {1 \over\lambda^{\sharp} \| \Sigma_{-1,-1}^{-1/2} W \zeta \|_2 } . $$
So now we have
$$ \Sigma_{-1,-1} c^0= { 1 \over\lambda^{\sharp}  \| \Sigma_{-1,-1}^{-1/2} \zeta \|_2^2} W \zeta $$
and hence
$$ \| \Sigma_{-1,-1} c^0\|_{\infty} = { \| w_{-S} \|_{\infty}  \over \lambda^{\sharp} \| \Sigma_{-1,-1}^{-1/2} \zeta \|_2^2} .$$
\hfill $\sqcup \mkern -12mu \sqcap$

{\bf Proof of Lemma \ref{non-sparse.lemma}.}
It holds that
$$ \Sigma_{-1,-1} (\gamma^0 - \gamma^{\sharp} ) =
\EE {\bf x}_{-1}^T {\bf x}_{-1} (\gamma^0 - \gamma^{\sharp}) =
\EE {\bf x}_{-1}^T ({\bf x}_{-S} {\rm A} {\bf x}_S ) \gamma_{-S}^0 = \Sigma_{-1,-1} c^0 .$$
So 
$$ \| \Sigma_{-1,-1} (\gamma^0 - \gamma^{\sharp} ) \|_{\infty} = {\| w_{-S} \|_{\infty}  \over \lambda^{\sharp}
\| \Sigma_{-1,-1}^{1/2} W\zeta \|_2^2 } ={\lambda^{\sharp} \| w_{-S} \|_{\infty} \over \lambda^{\sharp2}
\| \Sigma_{-1,-1}^{1/2}W  \zeta \|_2^2 } \le \lambda^{\sharp}
 . $$
Moreover
$$ \lambda^{\sharp} \| \gamma^{\sharp} \|_1 \le \sqrt {\rm s} \| \gamma^{\sharp} \|_2 \rightarrow 0 . $$
Thus $(\gamma^{\sharp}, \lambda^{\sharp})$ is an eligible pair. 
Furthermore
$$ \gamma^{0T} \Sigma_{-1, -1} \gamma^0 =
\gamma^{\sharp T} \Sigma_{-1,-1} \gamma^{\sharp} + c^{0T} \Sigma_{-1,-1} c^0 $$
$$= \gamma^{\sharp T} \Sigma_{-1,-1} \gamma^{\sharp} +
{ 1 \over \lambda^{\sharp2} \|  \Sigma_{-1,-1} ^{-1/2} W \zeta \|_2^2 } .  $$
So by Lemma \ref{allowedgamma.lemma} $\gamma^0$ is allowed. 
Finally, $\lambda^{\sharp} \| \gamma^0 \|_1 \ge \lambda^{\sharp} \| (W \gamma^0 )_{-S} \|_1 = 1 $
and in fact
$$ \gamma^{0T} \Sigma_{-1, -1} \gamma^0 -
\gamma^{\sharp T} \Sigma_{-1,-1} \gamma^{\sharp} ={1 \over  \lambda^{\sharp2} 
 \|  \Sigma_{-1,-1} ^{-1/2} W \zeta \|_2^2 } \gg 0  . $$

\hfill $\sqcup \mkern -12mu \sqcap$

\subsection{Proofs for Section \ref{unknownSigma.section}} 

{\bf Proof of Lemma \ref{slowLasso.lemma}.} 
We recall the notation $\hat \Sigma_{-1,-1} := X_{-1}^T X_{-1} / n $.
The event
$$\biggl \{  \| X_{-1}^T \epsilon^{\sharp} \|_{\infty} /n \le \lambda_{\varepsilon}^{\sharp} \biggr  \} \cap 
\biggl \{ \inf_{c: \ \lambda^{{\rm Lasso}} \| c \|_1 \le 4 \eta_n^2  , \ c^T \Sigma_{-1,-1} c = 1 } c^T \hat \Sigma_{-1,-1} c  \ge
{1 \over 2}  \biggr \}  $$
has probability converging to one so in the rest of the proof we may assume that
we are on this event. 
By the KKT conditions
$$ \hat \Sigma_{-1,-1} ( \hat \gamma - \gamma^0) = X_{-1}^T \varepsilon / n - \lambda^{{\rm Lasso}} \hat \zeta  $$
where $\hat \zeta \in \partial \| \hat \gamma \|_1 $, with $\partial \| c \|_1$ the sub-differential of the map
$c \mapsto \| c \|_1 $. 
Thus
$$ \hat \Sigma_{-1,-1} ( \hat \gamma - \gamma^{\sharp}) = X_{-1}^T \varepsilon^{\sharp} / n - 
\lambda^{{\rm Lasso}} \hat \zeta . $$
Therefore
\begin{eqnarray*}
( \hat \gamma - \gamma^{\sharp})^T \hat \Sigma_{-1,-1} ( \hat \gamma - \gamma^{\sharp}) &= &
( \hat \gamma - \gamma^{\sharp})^T X_{-1}^T \varepsilon^{\sharp} / n - \lambda^{{\rm Lasso}}( \hat \gamma - \gamma^{\sharp})^T\hat \zeta\\
&\le & \lambda_{\varepsilon}^{\sharp} \|  \hat \gamma - \gamma^{\sharp} \|_1+ \lambda^{{\rm Lasso}}
\| \gamma^{\sharp}\|_1 - \lambda^{{\rm Lasso}}\| \hat \gamma \|_1 \\
&\le & ( \lambda^{{\rm Lasso}} + \lambda_{\varepsilon}^{\sharp} )
\| \gamma^{\sharp}\|_1- (\lambda^{\rm Lasso} - \lambda_{\varepsilon}^{\sharp} )\| \hat \gamma \|_1 \\
&\le & ( \lambda^{{\rm Lasso}} + \lambda_{\varepsilon}^{\sharp} )
\| \gamma^{\sharp}\|_1\\
& \le& 3 \lambda^{\rm Lasso} \| \gamma^{\sharp} \|_1 /2 \le 3 \eta_n^2 /2 .
\end{eqnarray*}
It moreover follows from the above that
$$ \| \hat \gamma \|_1 \le \biggl ( {  \lambda^{{\rm Lasso}}  + \lambda_{\varepsilon}^{\sharp} \over \lambda^{{\rm Lasso}} - \lambda_{\varepsilon}^{\sharp}}
 \biggr )  \| \gamma^{\sharp}\|_1 $$
 and so 
 $$\lambda^{\rm Lasso} \| \hat \gamma  \|_1\le \lambda^{\rm Lasso} \biggl ( {  \lambda^{{\rm Lasso}} + \lambda_{\varepsilon}  \over 
\lambda^{{\rm Lasso}} - 
\lambda_{\varepsilon}^{\sharp}}\biggr ) \le 3 \eta_n^2 , $$
and also 
$$ \lambda^{\rm Lasso} \| \hat \gamma - \gamma^{\sharp} \|_1 \le \lambda^{\rm Lasso} \biggl ( {  2\lambda^{{\rm Lasso}}  \over 
\lambda^{{\rm Lasso}} - 
\lambda_{\varepsilon}^{\sharp}}
\biggr )  \| \gamma^{\sharp} \|_1 \le 4 \eta_n^2   . $$ 
If $\| \Sigma^{1/2} ( \hat \gamma - \gamma^{\sharp} )\|_2 \le 2 \eta_n$ we are done. 
Otherwise, if $\| \Sigma^{1/2} ( \hat \gamma - \gamma^{\sharp} )\|_2 \ge 2 \eta_n$
it holds that $4 \eta_n^2 \le 2 \eta_n \| \Sigma^{1/2} ( \hat \gamma - \gamma^{\sharp} )\|_2$.
But then
\begin{eqnarray*}
{1 \over 2}  ( \hat \gamma - \gamma^{\sharp})^T  \Sigma_{-1,-1} ( \hat \gamma - \gamma^{\sharp}) &\le&
( \hat \gamma - \gamma^{\sharp})^T  \hat \Sigma_{-1,-1} ( \hat \gamma - \gamma^{\sharp}) \\
&\le & (\lambda^{{\rm Lasso}} + \lambda_{\varepsilon}^{\sharp} )
\| \gamma^{\sharp}\|_1 \\
& \le  & 2 \lambda^{{\rm Lasso}} \| \gamma^{\sharp}\|_1  \le 2 \eta_n^2 .
\end{eqnarray*}
\hfill $\sqcup \mkern -12mu \sqcap$

{\bf Proof of Theorem \ref{noisy.theorem}.}
We rewrite
$$ \hat \beta_1 - \beta_1^0 = \hat \Theta_1 X^T \epsilon / n +
\underbrace{({\rm e}_1^T - \hat \Theta_1^T \hat \Sigma) ( \hat \beta- \beta^0) }_{\rm remainder} . $$
By the KKT conditions
$$ \| {\rm e}_1^T - \hat \Theta_1^T \hat \Sigma \|_{\infty} \le \lambda^{\rm Lasso}/(\| X_1 - X_{-1} \hat \gamma \|_2^2 / n +
 \lambda^{\rm Lasso} \| \hat  \gamma \|_1  ).$$
 But by Lemma \ref{slowLasso.lemma}
 $$\| X_1 - X_{-1} \hat \gamma \|_2^2 / n +
 \lambda^{\rm Lasso} \| \hat  \gamma \|_1 =
 \EE( {\bf x}_1 - {\bf x}_{-1} \gamma^{\sharp} )^2 + o(1) $$
 which stays away from zero.
 Moreover, by assumption, $\sqrt n \lambda^{\rm Lasso} \| \hat \beta - \beta_0 \|_1 = o_{\PP_{\beta^0}} (1)$
 uniformly in $\beta^0 \in {\cal B}$. 
 Thus, for the remainder we find
$$ |  ({\rm e}_1^T - \hat \Theta_1^T \hat \Sigma) ( \hat \beta- \beta^0) |\le
\| {\rm e}_1^T - \hat \Theta_1^T \hat \Sigma \|_{\infty} \| \hat \beta- \beta^0 \|_1 = o_{\PP_{\beta^0}} ( 1/ \sqrt n ) $$
uniformly in $\beta^0 \in {\cal B}$.
For the main term, we have after standardization
$$  { \hat \Theta_1^T X^T \epsilon /\sqrt  n \over
\sqrt {\hat \Theta_1^T \hat \Sigma \hat \Theta_1} } \sim {\cal N} (0,1) .$$
It further holds that  $\hat \Theta_1^T \hat \Sigma \hat \Theta_1 = \Theta_{1,1}^{\sharp} + o_{\PP} (1)$
by Lemma \ref{slowLasso.lemma}.
which stays away from zero.
Therefore, for  the standardized remainder term
$$ {\sqrt n  |  ({\rm e}_1^T - \hat \Theta_1^T \hat \Sigma) ( \hat \beta- \beta^0) | \over 
\sqrt {\hat \Theta_1^T \hat \Sigma \hat \Theta_1} } = o_{\PP_{\beta^0}} (1 )$$
uniformly in $\beta^0 \in {\cal B}$. 
The final result thus follows from Slutsky's Theorem.
\hfill $\sqcup \mkern -12mu \sqcap$

\section{Probability inequalities}\label{inequalities.section}

In this section we present some probability inequalities for products of Gaussians.
Such results are known (for example as Hanson-Wright inequalities for sub-Gaussians, see 
\cite{Rudelson2013}) and only presented here for completeness. 

 \begin{lemma} \label{UW.lemma}Let $U$ and $W$ be two independent ${\cal N} (0,1)$-distributed random variables.
   Then for all $L > 1$
   $$ \EE \exp \biggl [ {UW \over L} \biggr ] \le \exp \biggl [ {1 \over 2 L^2 - 2L } \biggr ] . $$
   \end{lemma}
   
   {\bf Proof.} We have for $L>1$
   \begin{eqnarray*}
    \EE \exp \biggl [ {UW \over L }\biggr ] &\le& \EE \exp \biggl [ {(U + W)^2 -(U-W)^2  \over 4L } \biggr ]\\
  & = &\EE\exp \biggl [ {(U+W)^2-2 \over 4L} \biggr ] \exp \biggl [ { 2- (U-W)^2  \over 4L} \biggr ] \\
 &= & \EE\exp \biggl [ {(U+W)^2-2 \over 4L} \biggr ] \EE  \exp \biggl [{ 2-(U-W)^2  \over 4L} \biggr ] \\
 & \le& \exp \biggl [ {1 \over 4L^2 - 4L}  \biggr ]  \EE  \exp \biggl [{ 2-(U-W)^2  \over 4L} \biggr ] .
  \end{eqnarray*}
 (see Lemma 1 and its proof in \cite{laurent2000adaptive}) or  Section  8.4 in \cite{vdG2016}).
  But
  \begin{eqnarray*}
   \EE  \exp \biggl [{ 2-(U-W)^2  \over 4L} \biggr ] &= &{1 \over \sqrt {2 \pi}}  \int \exp\biggl [ {1-v^2 \over 2L} \biggr ] \exp\biggl [ 
   - { v^2 \over 2 } \biggr ]  dv \\
   & = & \exp \biggl [ { 1 \over 2L} \biggr ] {1 \over \sqrt {2 \pi} } \int \exp \biggl [ - {1 \over 2} \biggl ( 1+ {1 \over L} 
   \biggr ) v^2 
   \biggr ]  dv \\
   & = & \exp \biggl [ { 1 \over 2L} \biggr ] \biggl ( 1+ {1 \over L} 
   \biggr )^{-1/2} \\ &= & \exp \biggl [ {1 \over 2L} - {1 \over 2} \log (1+ 1/L )  \biggr ] . 
   \end{eqnarray*}
   Since
   $$ \log (1+ 1/L) \ge  1/L - {{\ \atop 1} \over {2 \atop \ }} (1/L^2 )  $$
   we obtain
   \begin{eqnarray*}
   \EE  \exp \biggl [{ 2-(U+W)^2  \over 4L} \biggr ] &\le & \exp\biggl [ {1 \over 4 L^2 } \biggr ] \\
   &\le& \exp\biggl [ { 1 \over 4L^2 - 4L } \biggr ] . 
   \end{eqnarray*}
 It follows that
 $$  \EE \exp \biggl [ {UW  \over L }\biggr ]  \le \exp \biggl [ {2 \over 4L^2 - 4L}  \biggr ]  
 =  \exp \biggl [ {1 \over 2L^2 - 2L}  \biggr ] . $$
 \hfill $\sqcup \mkern -12mu \sqcap$ 
  
 \begin{lemma} \label{UWn.lemma} Let $ U= (U_1 , \ldots , U_n)^T $ and $W= (W_1 , \ldots , W_n)^T$ be two independent
standard Gaussian $n$-dimensional random vectors.
Then for all $t >0$
$$ \PP \biggl ( U^T W /n \ge  \sqrt {2t/n} +  t/n \biggr ) \le \exp [-t] . $$
 \end{lemma}
 
 {\bf Proof.} By Lemma \ref{UW.lemma} and using the independence
 $$ \EE \exp \biggl [ { U^T W   \over L} \biggr ] \le
 \exp \biggl [ {n \over 2L^2 - 2L}  \biggr ]. 
 $$
 This gives for all $t>0$
 $$ \PP \biggl ( U^T W \ge  \sqrt {2nt} +  t \biggr ) \le \exp [-t]  $$
 (see e.g.\  Lemma 8.3 in \cite{vdG2016}).
 \hfill $\sqcup \mkern -12mu \sqcap$ 
   
   \begin{lemma} \label{lambda*.lemma}Let $\{ (U_i, V_i) \}_{i=1}^n $ be i.i.d.\ two-dimensional Gaussians with mean zero.
    Suppose ${\rm var} (U_1) =1$. Define $\lambda^{\sharp} := \EE U_1V_1 $ and
   $\sigma^{\sharp 2}= \EE V_1^2$.  Then for all $t >0$
   $$\PP \biggl ( U^T V / n \ge \lambda^{\sharp} + ( \sqrt 2  \sigma^{\sharp} + 2\lambda^{\sharp}) \sqrt {t/n} +
   ( \sigma^{\sharp} + 2 \lambda^{\sharp}) t/n \biggr ) \le 2 \exp[-t] . $$

   \end{lemma}
   
   {\bf Proof.} For all $i$ the projection of $V_i$ on $U_i$ is $[\EE U_iV_i / {\rm var } (U_i) ] U_i = \lambda^{\sharp} U_i$. Hence
   we may write for all $i$
   $$ V_i = \lambda^{\sharp} U_i + W_i , $$
   where $W_i$ is a zero-mean Gaussian random variable independent of $U_i$. It follows
   that
   $$ U^T V/n  = \lambda^{\sharp} \| U \|_2^2 /n + U^T W/n . $$
   Since ${\rm var} ( W_i ) \le \sigma^{\sharp2} $ for all $i$ we see from Lemma \ref{UWn.lemma}
   that
   $$ \PP \biggl ( U^T W / n \ge  \sqrt 2\sigma^{\sharp} \sqrt { t/n} + \sigma^{\sharp} t/n \biggr ) \le \exp [-t] . $$
   Moreover (see Lemma in \cite{laurent2000adaptive}, also given in \cite{vdG2016} as Lemma 8.6)
   $$ \PP\biggl  ( \| U \|_2^2/n  -1 \ge 2 \sqrt {t/n} + 2 t/n\biggr  ) \le \exp[-t] . $$
   Thus
   $$\PP \biggl ( U^T V / n \ge \lambda^{\sharp} + (  \sqrt 2 \sigma^{\sharp} + 2 \lambda^{\sharp}) \sqrt {t/n} +
   ( \sigma^{\sharp} + 2 \lambda^{\sharp}) t/n \biggr ) \le 2 \exp[-t] . $$
   \hfill $\sqcup \mkern -12mu \sqcap$

\bibliographystyle{plainnat}
\bibliography{reference}
\end{document}